\RequirePackage{fix-cm}
\documentclass{svjour3}                     
\smartqed  
\usepackage{graphicx}
\RequirePackage[OT1]{fontenc}
\RequirePackage[numbers]{natbib}
\RequirePackage[colorlinks,citecolor=blue,urlcolor=blue]{hyperref}
\usepackage{graphicx}

  \usepackage[dvipsnames]{xcolor}
\usepackage[utf8]{inputenc}
\usepackage[english]{babel}
\usepackage[T1]{fontenc}
\usepackage{amssymb, amsfonts,  amsmath}
\usepackage{dsfont}
\usepackage{bbm}
\usepackage{tikz}
\usepackage{longtable}
  \usepackage{graphicx}
\usepackage{hyperref}
\usepackage{xspace}
\usepackage[plain,noend]{algorithm2e}
\usepackage{placeins}
  \colorlet{mred}{red!80!black}
\definecolor{mblue}{RGB}{94,135,173}
\usepackage{placeins}
\usepackage{multirow}
\usepackage{array}
\newcounter{assumptA}

\newcounter{assumptM}

\spdefaulttheorem{nota}{Notations}{\bf}{\it}
\spdefaulttheorem{assumpt}{Assumptions}{\bf}{\it}

\newcommand{\dd}{d}
\newcommand{\g}{g}
\newcommand{\ii}{i}
\newcommand{\jj}{j}
\newcommand{\kk}{k}
\newcommand{\el}{\ell}
\newcommand{\m}{m}
\newcommand{\n}{n}
\newcommand{\vrai}{\star}
\newcommand{\al}{\alpha}
\newcommand{\als}{\al^\vrai}
\newcommand{\alkl}{\al_{\kk\el}}
\newcommand{\alskl}{\als_{\kk\el}}
\newcommand{\bal}{\boldsymbol{\al}}
\newcommand{\bals}{\boldsymbol{\al}^\vrai}

\newcommand{\thetaa}{\theta}

\newcommand{\btheta}{\boldsymbol{\thetaa}}

\newcommand{\bthetas}{\boldsymbol{\thetaa}^\vrai}

\newcommand{\pii}{\pi}
\newcommand{\bpi}{\boldsymbol{\pii}}
\newcommand{\pik}{\pii_\kk}
\newcommand{\bpis}{\boldsymbol{\pii}^\vrai}
\newcommand{\rhoo}{\rho}
\newcommand{\rhol}{\rhoo_\el}
\newcommand{\brho}{\boldsymbol{\rhoo}}
\newcommand{\brhos}{\boldsymbol{\rhoo}^\vrai}

\newcommand{\w}{w}
\newcommand{\bw}{\mathbf{\w}}
\newcommand{\bfW}{\mathbf{W}}
\newcommand{\bws}{\bfw^\vrai}
\newcommand{\bWs}{\bfW^\vrai}
\newcommand{\wjl}{\w_{\jj\el}}
\newcommand{\Wjl}{W_{\jj\el}}
\newcommand{\Wsjl}{\Wjl^\vrai}
\newcommand{\x}{x}
\newcommand{\bx}{\mathbf{\x}}
\newcommand{\xij}{\x_{\ii\jj}}
\newcommand{\X}{X}
\newcommand{\Xij}{\X_{\ii\jj}}
\newcommand{\z}{z}
\newcommand{\bfz}{\mathbf{\z}}
\newcommand{\bfZ}{\mathbf{Z}}
\newcommand{\bz}{\bfz}
\newcommand{\bZ}{\mathbf{Z}}
\newcommand{\bW}{\mathbf{W}}
\newcommand{\bfw}{\bw}

\newcommand{\zik}{\z_{\ii\kk}}
\newcommand{\Zik}{Z_{\ii\kk}}
\newcommand{\bzs}{\bfz^\vrai}
\newcommand{\bZs}{\bfZ^\vrai}

\newcommand{\prob}{p}
\newcommand{\Prob}{\mathbb{P}}

\newcommand{\EspC}[2]{\mathbb{E}_{#1}\left[#2\right]}

\newcommand{\dens}{\varphi}

\newcommand{\mcZ}{\mathcal{Z}}
\newcommand{\mcW}{\mathcal{W}}

\newcommand{\sommecolonne}{+}
\newcommand{\zsumk}{\z_{\sommecolonne\kk}}
\newcommand{\wsuml}{\w_{\sommecolonne\el}}

\newcommand{\zs}{\z^{\vrai}}
\newcommand{\ws}{\w^{\vrai}}
\newcommand{\Zs}{Z^{\vrai}}
\newcommand{\Ws}{W^{\vrai}}
\newcommand{\zssumk}{\zs_{\sommecolonne\kk}}
\newcommand{\wssuml}{\ws_{\sommecolonne\el}}
\newcommand{\Zssumk}{\Zs_{\sommecolonne\kk}}
\newcommand{\Wssuml}{\Ws_{\sommecolonne\el}}

\newcommand{\kp}{\kk'}

\newcommand{\transpose}{{^{T}}}

\newcommand{\lp}{\el'}

\newcommand{\G}{G}

\newcommand{\norminf}[1]{\left\|#1\right\|_{\infty}}

\newcommand{\y}{y}

\newcommand{\tauu}{\tau}
\newcommand{\btau}{\boldsymbol{\tauu}}
\newcommand{\btaus}{\btau^{\vrai}}
\newcommand{\bxis}{\bxi^{\vrai}}
\newcommand{\LG}{}
\newcommand{\tauk}{\tauu_{\kk}}
\newcommand{\somme}{+}
\newcommand{\Xiplus}{\X_{\ii\somme}}
\newcommand{\zplusk}{\z_{\somme\kk}}
\newcommand{\wplusl}{\w_{\somme\el}}
\newcommand{\Bin}{\mathcal{B}in}
\newcommand{\xii}{\xi}
\newcommand{\bxi}{\boldsymbol{\xii}}
\newcommand{\Tid}{\overline{\Xidot}}
\newcommand{\Tund}{\overline{\X_{(1)\ligne}}}
\newcommand{\Tnun}{\overline{\X_{\ligne(1)}}}
\newcommand{\Tnd}{\overline{\X_{(\n)\ligne}}}
\newcommand{\Tndbis}{\overline{\X_{\ligne(\dd)}}}
\newcommand{\Tipard}{\overline{\X_{(\ii)\ligne}}}
\newcommand{\ligne}{\centerdot}
\newcommand{\seuil}{S}
\newcommand{\Sg}{\seuil_{\g}}
\newcommand{\Sm}{\seuil_{\m}}
\newcommand{\bzLG}{\widehat{\bZ}^{\LG}}
\newcommand{\bwLG}{\widehat{\bW}^{\LG}}
\newcommand{\bX}{\mathbf{\X}}
\newcommand{\Xidot}{\X_{\ii\ligne}}
\newcommand{\Tiparund}{\overline{\X_{(\ii-1)\ligne}}}
\newcommand{\gLG}{\widehat{\g}^{\LG}}
\newcommand{\Gi}{\G_{\ii}}
\newcommand{\mLG}{\widehat{\m}^{\LG}}
\newcommand{\pikLG}{\widehat{\pik}^{\LG}}
\newcommand{\rholLG}{\widehat{\rhol}^{\LG}}
\newcommand{\balLG}{\widehat{\bal}^{\LG}}
\newcommand{\bthetaLG}{\widehat{\btheta}}
\newcommand{\permbthetaLG}{\widehat{\btheta}^{\perms_{\clZ},\permt_{\clW}}}
\newcommand{\zLG}{\widehat{\z}^{\LG}}
\newcommand{\wLG}{\widehat{\w}^{\LG}}
\newcommand{\alklLG}{\widehat{\alkl}^{\LG}}
\newcommand{\zikLG}{\widehat{\zik}^{\LG}}
\newcommand{\wjlLG}{\widehat{\wjl}^{\LG}}
\newcommand{\zsumkLG}{\widehat{\zsumk}^{\LG}}
\newcommand{\wsumlLG}{\widehat{\wsuml}^{\LG}}
\newcommand{\pimin}{\pii_{\min}^{\vrai}}
\newcommand{\rhomin}{\rhoo_{\min}^{\vrai}}
\newcommand{\gs}{\g^{\vrai}}
\newcommand{\piks}{\pik^{\vrai}}
\newcommand{\ms}{\m^{\vrai}}
\newcommand{\rhols}{\rhol^{\vrai}}
\newcommand{\deltaa}{\delta}
\newcommand{\deltapi}{\deltaa_{\bpi^\vrai}}
\newcommand{\deltarho}{\deltaa_{\brho^\vrai}}
\newcommand{\taukp}{\tauu_{\kp}}
\newcommand{\xil}{\xii_{\el}}
\newcommand{\xilp}{\xii_{\lp}}
\newcommand{\Z}{Z}
\newcommand{\W}{W}
\newcommand{\clW}{\mathcal{\W}}
\newcommand{\clZ}{\mathcal{\Z}}
\newcommand{\egalZ}{\equiv_{\clZ}}
\newcommand{\cegalZ}{\not\equiv_{\clZ}}
\newcommand{\egalW}{\equiv_{\clW}}
\newcommand{\cegalW}{\not\equiv_{\clW}}
\newcommand{\dHinf}{d^{\infty}}
\newcommand{\by}{\mathbf{\y}}

\newcommand{\perms}{s}
\newcommand{\permt}{t}
\newcommand{\Sgnd}{\Sg^{\n,\dd}}
\newcommand{\Smnd}{\Sm^{\n,\dd}}
\newcommand{\deltapind}{\deltapi^{\n,\dd}}
\newcommand{\deltarhond}{\deltarho^{\n,\dd}}
\newcommand{\piminnd}{\pii_{\min}^{\vrai\,\n,\dd}}
\newcommand{\rhominnd}{\rhoo_{\min}^{\vrai\,\n,\dd}}
\newcommand{\gsnd}{{\gs}^{\n,\dd}}
\newcommand{\msnd}{{\ms}^{\n,\dd}}
\newcommand{\Sun}{\seuil_{1}}
\newcommand{\Sdeux}{\seuil_{2}^{\n,\dd}}
\newcommand{\Strois}{\seuil_{3}^{\n,\dd}}
\newcommand{\Squatre}{\seuil_{4}^{\n,\dd}}
\newcommand{\Ags}{A_{\gs}}
\newcommand{\D}{D}
\newcommand{\zsik}{\z_{\ii,\kk}^{\vrai}}
\newcommand{\ASg}{A_{\Sg}}

\newcommand{\Aid}{A_{id}}
\newcommand{\ip}{\ii'}
\newcommand{\Tipd}{\overline{\Xipdot}}
\newcommand{\Xipdot}{\X_{\ip\ligne}}
\newcommand{\Zsik}{\Z_{\ii,\kk}^{\vrai}}
\newcommand{\bpiLG}{\widehat{\bpi}^{\LG}}
\newcommand{\brhoLG}{\widehat{\brho}^{\LG}}
\newcommand{\alkls}{\alkl^{\vrai}}
\newcommand{\wsjl}{\w_{\jj,\el}^{\vrai}}
\newcommand{\alkltilde}{\widetilde{\alkl}}

\newcommand{\N}{N}
\newcommand{\tauks}{\tauk^{\vrai}}
\newcommand{\taukps}{\taukp^{\vrai}}
\newcommand{\xils}{\xil^{\vrai}}
\newcommand{\xilps}{\xilp^{\vrai}}

\newcommand{\Rev}[1]{#1}

\allowdisplaybreaks

\begin{document}

\title{Fast and Consistent Algorithm for the Latent Block Model\thanks{Thanks to St\'{e}phane Robin and the reviewers for them suggestions. This work has been partially supported by the LabEx PERSYVAL-Lab (ANR-11-LABX-0025-01) funded by the French program Investissement d’avenir and partially supported by MIAI@Grenoble Alpes (ANR-19-P3IA-0003). All the computations presented in this paper were performed using the GRICAD infrastructure (\url{https://gricad.univ-grenoble-alpes.fr}), which is supported by Grenoble research communities.
}
}


\author{Vincent Brault         \and
        Antoine Channarond 
}


\institute{Vincent Brault \at
              Univ. Grenoble Alpes, Inria, CNRS, Grenoble INP\footnote{Institute of Engineering Univ. Grenoble Alpes}, LJK, 38000 Grenoble, France \\
              \email{vincent.brault@univ-grenoble-alpes.fr}           
           \and
           Antoine Channarond \at
              UMR6085 CNRS, Laboratoire de Mathématiques Raphaël Salem, Université de Rouen Normandie, 76800 Saint-\'Etienne-du-Rouvray, France\\
              \email{antoine.channarond@univ-rouen.fr}
}

\date{Received: date / Accepted: date}

\maketitle

\begin{abstract}
The latent block model is used to simultaneously rank the rows and columns of a matrix to reveal a block structure. The algorithms used for estimation are often time consuming. However, recent work shows that \Rev{the log-likelihood ratios are equivalent under the complete and observed (with unknown labels) models and the groups posterior distribution to converge as the size of the data increases to a Dirac mass located at the actual groups configuration}. Based on these observations, the algorithm \textit{Largest Gaps} is proposed in this paper to perform clustering using only the marginals of the matrix, when the number of blocks is very small with respect to the size of the whole matrix in the case of binary data. In addition, a model selection method is incorporated with a proof of its consistency. \Rev{Thus, this paper shows that studying simplistic configurations (few blocks compared to the size of the matrix or very contrasting blocks) with complex algorithms is useless since the marginals already give very good parameter and classification estimates.}

\keywords{Latent Block Model \and Largest Gaps Algorithm \and Model Selection \and Data analysis}
\subclass{\\62H30 \and 62-07 \and 60K35}
\end{abstract}

\section{Introduction}

Block clustering methods aim at clustering rows and columns of a matrix simultaneously to form homogeneous blocks. There are a lot of applications of this method: genomics \citep{hedenfalk2001gene,jagalur2007analyzing}, recommender systems \citep{bennett2007netflix,shan2008bayesian}, archeology \citep{govaert1983classification}, sociology \citep{hartigan1975Clustering,keribin2014estimation,wyse2010block} or network \citep{barbillon2017stochastic,bouveyron2018stochastic} for example. Among the methods proposed to solve this question, the Latent Block Model (LBM) \citep{govaert2003clustering} provides a chessboard structure induced by the classification of the rows and the classification of the columns. In this model, it is assumed that a sample of $\n$ individuals is collected, which contains the observation of $\dd$ binary variables of the same nature. Saying that the binary variables are of the same nature means that it is possible to encode them in the same (and natural) way. This assumption is needed to ensure that decomposing the dataset in a block structure makes sense \Rev{and \cite{laclau2017co} shows the equivalence between searching for a classification in the LBM and solving an optimal transport problem}. In the case where the rows and columns represent the same individuals, we \Rev{generally} speak instead of \textit{Stochastic Block Model} \citep{barbillon2017stochastic,bouveyron2018stochastic,tabouy2020variational}. \Rev{However, in the case where the matrices are not symmetric, \cite{keribin2021cluster} shows that it may be more interesting to use the LBM model rather than the SBM model.} Although this model is very similar, we will not discuss it in this article. 

Given the number of blocks and in order to estimate the parameters, \cite{govaert2003clustering} suggest using a variational algorithm, \cite{keribin2012model} propose an adaptation of the Stochastic Expectation Maximisation introduced by \cite{celeux1995stochastic} in the mixture case, \cite{keribin2014estimation} studied a Bayesian version of these two algorithms and \cite{wyse2010block} propose a Bayesian algorithm including the estimation of the number of blocks. However, each of these iterative algorithms has a complexity at least $\mathcal{O}\left(\n\dd\N_{Block}\N_{Iter}\right)$ where $\N_{Block}$ is  the number of blocks and $\N_{Iter}$ is the number of iterations necessary for the convergence of the algorithm. Moreover, the procedures are often associated with a model selection criterion requiring the computation of maximum likelihood estimators for all combinations of the expected number of blocks.
However, the theoretical results obtained show that the distributions of the estimators (see \cite{brault2020consistency,celisse2012consistency,mariadassou2012convergence}) are asymptotically trivial: \Rev{the log-likelihood ratios are equivalent under the complete and observed (with unknown labels) models and the groups posterior distribution to converge as the size of the data increases to a Dirac mass located at the actual groups configuration}.

In this article, we propose an adaptation of the \textit{Largest Gaps} (LG) algorithm introduced by \cite{channarond2011Classification} in the \textit{Stochastic Block Model} with a low complexity $\mathcal{O}(\n\dd)$ (Section~\ref{sec:algo}) and \Rev{ present the conditions to obtain asymptotically a good estimation (Section~\ref{sec:assumptions:model}). We then} prove that it provides a consistent procedure for all inference tasks inherent in LBM: unsupervised classification and estimation of the parameters and a selection procedure for the number of blocks (the last unknown theoretical point) is also proposed and shown to be consistent (Section~\ref{sec:consistency}). For ease of reading, the proofs of the results are postponed to the appendices \Rev{(Sections~\ref{Appendix:thm:concentration}, \ref{cor:Appendix:consistency} and~\ref{cor:Appendix:consistency:spare})}. These theoretical results are also illustrated on simulated data (Section~\ref{sec:simulation}).

By proving the consistency of the \textit{LG} algorithm, some features of the asymptotic regime of the LBM will be highlighted, in particular, the concentration of some marginal distributions of the model. The secondary objective of the paper is to discuss as a conclusion (Section~\ref{sec:conclusion}) the consequence of this when fitting the LBM to large matrices; in particular, this leads to take a step back regarding the relevance of the model and the estimates when the number of blocks is small with respect to the number of cells of the matrix.

\section{Notations and model}

The binary Latent Block Model (LBM) is as follows.
Let $\bx=\left(\xij\right)_{\ii=1, \ldots, \n; \jj=1,\ldots, \dd}$ be the data matrix
where $\xij\in\{0,1\}$; \Rev{observation of a random variable $\bX$}. It is assumed that there exists a partition into $\g$ row clusters $\bz=\left(\zik\right)_{\ii=1, \ldots, \n; \kk=1,\ldots, \g}$ and a partition into $\m$ column clusters
$\bw=\left(\wjl\right)_{\jj=1, \ldots, \dd; \el=1,\ldots, \m}$.
The $\zik$s (resp. $\wjl$s ) are binary indicators of row $\ii$ (resp. column $\jj$) belonging to row cluster $\kk$ (resp. column cluster $\el$), such that
the random variables $\Xij$ are independent conditionally on $\bz$ and $\bw$ with parametric density $\dens(\cdot;\alkl)^{\zik \wjl}$, where $\alkl$ is the parameter of the conditional density of the data given $\Zik=1$ and $\Wjl=1$. \Rev{As the binary case is studied, the density is defined for each $\x\in\{0,1\}$ and $\alpha\in[0;1]$ by}
\[\Rev{\dens (\x;\al)=\al^{\x}\left(1-\al\right)^{1-\x}.}\]
Thus, the density of $\bX$ conditionally on $\bz$ and $\bw$ is defined for each $\bx$ by
\begin{eqnarray*}
f(\bx|\bz,\bw;\bal)&=&\prod_{\ii=1}^\n\prod_{\jj=1}^\dd\prod_{\kk=1}^\g\prod_{\el=1}^\m \dens(\xij;\alkl)^{\zik \wjl}\\
&=:&\prod_{\ii,\jj,\kk,\el} \dens(\xij;\alkl)^{\zik \wjl}
\end{eqnarray*}
where $\bal=\left(\alkl\right)_{\kk=1, \ldots, \g; \el=1,\ldots, \m}$.

Moreover, it is assumed that the row and column labels $\bz$ and $\bw$ are the observations of two independent variables $\bZ$ and $\bW$: $p(\bz,\bw;\bpi,\brho)=p(\bz;\bpi)p(\bw;\brho)$ with
$\prob(\bz;\bpi)=\prod_{\ii,\kk} \pik^{\zik}=\prod_{\scriptscriptstyle k} \pik^{\zplusk}$ and $\prob(\bw;\brho)=\prod_{\jj,\el} \rhol^{\wjl}=\prod_{\scriptscriptstyle \ell} \rho_{\ell}^{\wplusl}$, where $(\pik= \Prob(\Zik=1),k=1,\ldots,g)$ and $(\rhol=\Prob(\Wjl=1),\ell=1,\ldots,m)$ are the mixing proportions and  $\zplusk=\sum_{\ii=1}^\n{\zik}$ (resp. $\wplusl=\sum_{\jj=1}^\dd{\wjl}$) represents the number of rows (resp. columns) in the class $\kk$ (resp. $\ell$). Hence, the density of $\bX$ is defined for every $\bx$ by
\begin{equation*}
f(\bx;\btheta)=\sum_{(\bz,\bw) \in \mcZ\times\mcW}
p(\bz;\bpi)p(\bw;\brho)f(\bx|\bz,\bw;\bal),
\end{equation*}
where $\mcZ$ and $\cal W$ denoting the sets of all possible row labels $\bz$ and column labels $\bw$,
and $\btheta=(\bpi,\brho,\bal)$, with $\bpi=(\pi_1,\ldots,\pi_g)$ and
$\brho=(\rho_1,\ldots,\rho_m)$.
The density of $\bX$ can be written for every $\bx$ as
\begin{eqnarray*}
    f(\bx;\btheta)&=& \sum_{(\bz,\bw) \in \mcZ\times\mcW}
\prod_{\scriptscriptstyle k} \pik^{\zplusk} \prod_{\scriptscriptstyle \ell} \rho_{\ell}^{\wplusl}  \!\prod_{\scriptscriptstyle \ii,\jj,\kk,\el}\dens (\xij;\alkl)^{\zik\wjl}.\label{densite}\\\nonumber
\end{eqnarray*}

To estimate both the classification and the parameters, many algorithms exist (for example~\cite{govaert2003clustering}, \cite{keribin2014estimation} or \cite{wyse2010block}) but each of these algorithms has a complexity larger than $\mathcal{O}\left(\n\dd\g\m\N_{Iter}\right)$ where $\N_{Iter}$ is the number of iterations of the algorithm. This makes their use on large matrices difficult.

In the Stochastic Block Model (SBM), rows and columns are associated with the same individuals, which allows to represent a graph, whereas LBM allows to represent \Rev{bipartite graphs}. \cite{channarond2011Classification} suggested a fast algorithm, called $LG$ for \textit{Lagest Gaps}, based on a marginal of the matrix $\bx$, the degrees.

\section{Algorithm \textit{Largest Gaps}}\label{sec:algo}

Before the introduction of the classification algorithm \textit{Largest Gaps} ($LG$), let us recall the main idea, inspired by \cite{channarond2011Classification}.

\subsection{Main ideas} \label{subsec:concept}

Conditionally on $\Zik=1$, meaning that row $\ii$ is in class~$\kk$, the probability that the variable $\Xij$ equals $1$ for any $j\in\{1,\ldots,d\}$ is
\begin{eqnarray}
\Prob\left(\Xij=1|\Zik=1\right)
&=&\sum_{\el=1}^{\m}\Prob\left(\Xij=1|\Zik=1,\Wjl=1\right)\Prob\left(\Wjl=1|\Zik=1\right)\nonumber\\
&=&\sum_{\el=1}^{\m}\Prob\left(\Xij=1|\Zik=1,\Wjl=1\right)\Prob\left(\Wjl=1\right)\nonumber\\
&=&\sum_{\el=1}^{\m}\alkl\rhol\Rev{=:}\tauk.\label{Eq:Concept}
\end{eqnarray}
In particular, conditionally on $\Zik=1$, variables of row $\ii$ are independent and identically distributed Bernoulli variables with parameter $\tauk$ and the sum of the cells of any row $\ii$, denoted by $\Xiplus$, is hence a binomial distributed variable $\Bin\left(\dd,\tauk\right)$.

As a consequence of the subgaussian concentration property of binomial distributions, variables $\Tid=\frac{\Xiplus}{\dd}$ fastly concentrates around the mean associated with its own class when $\dd$ tend to infinity. The point is here: if moreover these means $\tauu_1,\ldots,\tauu_\g$ are pairwise distinct, the set of the variables $\{\Tid;\ii=1,\dots,\n\}$ asymptotically splits into clusters, separated by large gaps, and which exactly correspond to the clusters of the model, for $\dd$ large enough. In the whole paper, $\tauu_1,\ldots,\tauu_\g$ will be assumed to be pairwise distinct (Assumption \eqref{hyp:identifiability}, see discussion in Section \ref{sec:assumptions:model}).

The middle right picture of Figure~\ref{Fig:Improvement} shows the histogram of the variable set $\{\Tid;\ii=1,\dots,\n\}$ for a matrix simulated under LBM with five clusters. The five clusters of rows can be seen, as well as the four large gaps which separate them. The middle left picture of Figure~\ref{Fig:Improvement} is a representation of the vector $(\Tid$ sorted in ascending order and the bottom left picture of Figure~\ref{Fig:Improvement}, the size of the gaps between \Rev{two consecutive sorted values}.

To classify the rows, the idea is to identify the gaps between the clusters, which are expected to be asymptotically larger. Indeed, the internal gaps (those between two rows of the same cluster) vanish when $\dd$ tends to infinity due to concentration, while the external gaps (those between two rows of distinct clusters) do not, since the means are assumed to be distinct. Symmetrically the same holds for the columns.


There are several strategies to identify the large gaps. In their article, \cite{channarond2011Classification} assume that the number $Q$ of clusters (the same for the rows and the columns) is known and partition the population into $Q$ clusters by finding the $Q-1$ largest gaps. In order to choose $Q$, a model selection procedure can be firstly done separately, before the classification. The strategy chosen in this article consists in thresholding the gaps, to distinguish outer gaps from inner ones. It advantageously yields both the clusters and the numbers of clusters in only one pass.

The choice of thresholds is critical. Let's comment on the beautiful configuration of Figure~\ref{Fig:Improvement}, where the asymptotic regime has been reached, to present the key issues. If the thresholds are too large (greater than 0.06 on that example), the thresholding step will select only some of the four outer deviations and will not distinguish the five clusters from each other. Conversely, if they are too small (less than 0.01), the thresholding step will select all four outer gaps, but also some undesirable inner gaps, and the algorithm will wrongly separate some clusters.

\subsection{Algorithm}\label{sec:algorithm}



\begin{figure*}[!h]
    \begin{center}
        \begin{tabular}{cc}
            \begin{minipage}[c]{0.45\textwidth}
                \includegraphics[width=\linewidth]{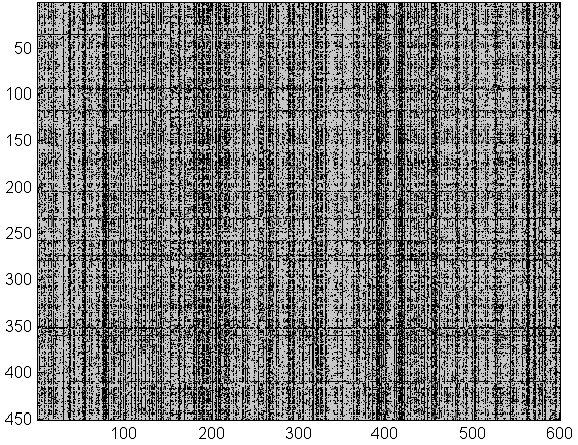}
            \end{minipage}&
            \begin{minipage}[c]{0.45\textwidth}
                \includegraphics[width=\linewidth]{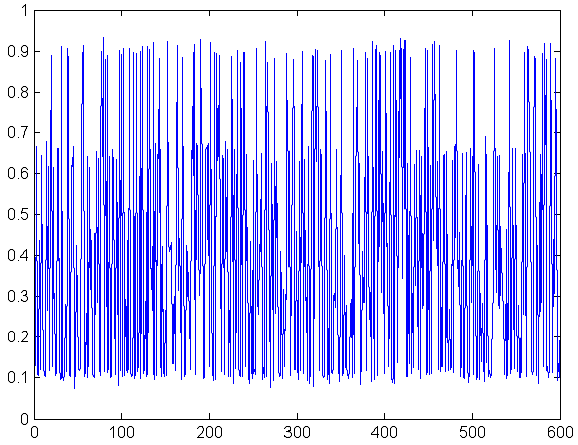}
            \end{minipage}\\
            \begin{minipage}[c]{0.45\textwidth}
                \includegraphics[width=\linewidth]{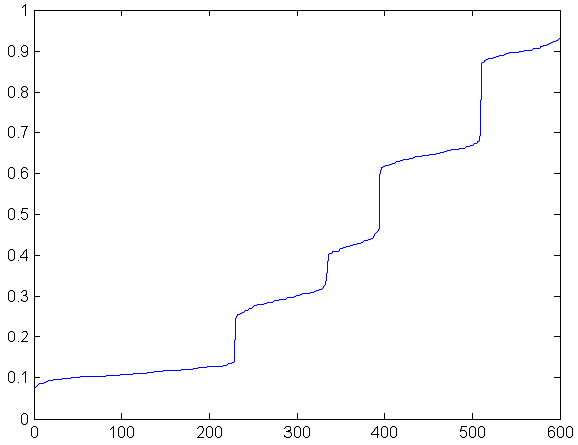}
            \end{minipage}&
            \begin{minipage}[c]{0.45\textwidth}
                \includegraphics[width=\linewidth]{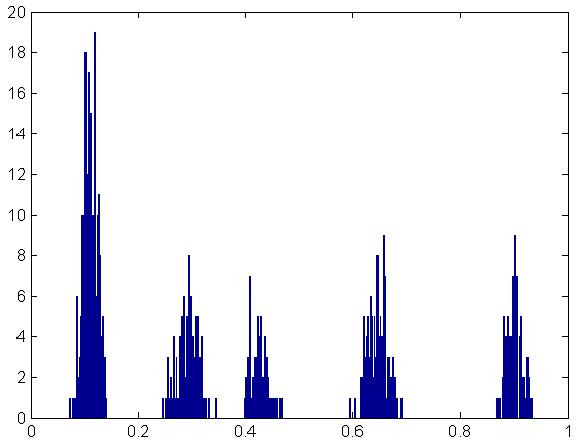}
            \end{minipage}\\
            \begin{minipage}[c]{0.45\textwidth}
                \includegraphics[width=\linewidth]{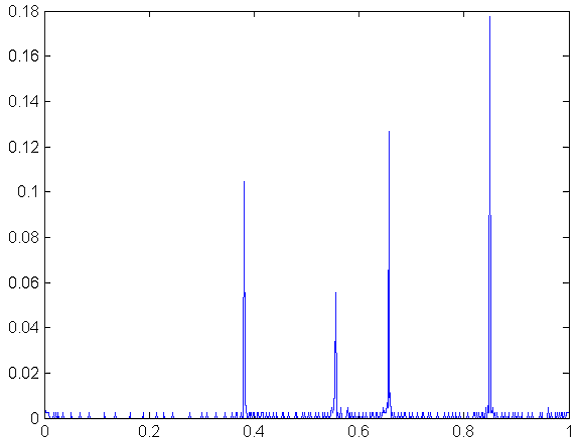}
            \end{minipage}&
            \begin{minipage}[c]{0.45\textwidth}
                \includegraphics[width=\linewidth]{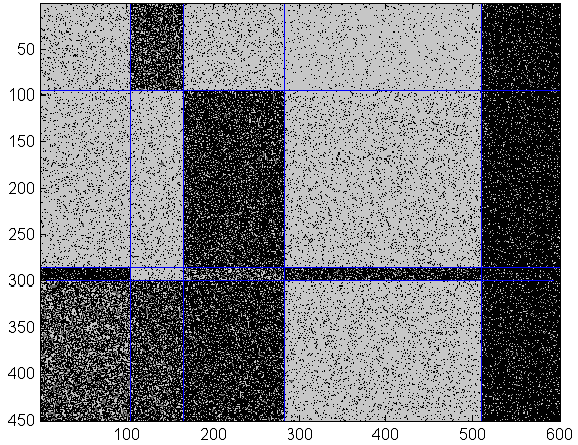}
            \end{minipage}\\
        \end{tabular}
    \end{center}
    \caption{\label{Fig:Improvement} Top-left: Initial matrix. Top-right: Example of a vector $\left(\Tnun,\ldots,\Tndbis\right)$. Middle-left: representation of the vector $\left(\Tnun,\ldots,\Tndbis\right)$ sorted in increasing order. Middle-right: Histograms of $\left(\Tnun,\ldots,\Tndbis\right)$. Bottom-left: representation of the vector of gaps $\left(G_2,\ldots,G_{\dd}\right)$ where for all $\jj\in\{2,\ldots,\dd\}$, $G_\jj=\overline{\X_{\ligne(\jj)}}-\overline{\X_{\ligne(\jj-1)}}$. Bottom-right: reorganized matrix.}
\end{figure*}

The algorithm \textit{Largest Gaps}\footnote{\Rev{An implementation in the language \texttt{R} is available on the following Gitlab: \url{https://gricad-gitlab.univ-grenoble-alpes.fr/braultv/largest-gaps}}} is given in \Rev{Algorithm~\ref{algo:LG}} and an illustration is provided in Figure~\ref{Fig:Improvement}. The principle is to calculate the means of the values of the cells of each row and each column, to order them and to calculate the differences between two consecutive values. Once this step is done, the clusters are formed assuming that each difference greater than the threshold implies a change of cluster. In the sequel, the estimators provided by the algorithm are denoted by $\bzLG$, $\bwLG$ and $\bthetaLG$.\\

\begin{algorithm}[htbp!]
  {
  \SetSideCommentLeft
  \DontPrintSemicolon
  \KwIn{data matrix $\bx$, threshold for row $\Sg$ and for column $\Sm$.}
  \BlankLine
  \BlankLine
  \tcp{\textcolor{mred}{Computation of gaps}}
  \For{$\ii\in\{1,\ldots,\n\}$} {
   Computation of $\Tid=\frac{x_{\ii\sommecolonne }}{\dd}$.
   }\tcp*[r]{\textcolor{mblue}{$\mathcal{O}(\n\dd)$}}
  \BlankLine
  Ascending sort of $\left(\Tund ,\ldots,\Tnd \right)$.\tcp*[r]{\textcolor{mblue}{$\mathcal{O}(\n\log\n)$}}
  \BlankLine
  \For{$\ii\in\{2,\ldots,\n\}$} {
   Computation of the gaps $\Gi=\Tipard-\Tiparund$.
   }\tcp*[r]{\textcolor{mblue}{$\mathcal{O}(n)$}}
  \BlankLine
  \tcp{\textcolor{mred}{Computation of $\gLG$}}
  Selection of $\ii_1<\ldots<\ii_{\gLG-1}$ such that $(G_{\ii_1},\ldots,G_{\ii_{\gLG-1}})$ are every greater than $\Sg$.\tcp*[r]{\textcolor{mblue}{$\mathcal{O}(\n)$}}
  \BlankLine
  \tcp{\textcolor{mred}{Computation of $\bzLG$}}
  \For{$\ii\in\{(1),\ldots,(\n)\}$} {
  Definition of $\zLG_{(i)k}=1$ if and only if $(\ii_{\kk-1})<(\ii)\leq(\ii_\kk)$ with $\ii_0=0$ and $\ii_{\gLG}=\n$.
  }\tcp*[r]{\textcolor{mblue}{$\mathcal{O}(n)$}}
  \BlankLine
  \tcp{\textcolor{mred}{Computation of $\mLG$ and $\bwLG$}}
  Do the same on the columns.\tcp*[r]{\textcolor{mblue}{$\mathcal{O}(\max(\dd\n,\dd\log\dd))$}}
  \BlankLine
  \tcp{\textcolor{mred}{Computation of $\bthetaLG$}}
  \BlankLine
  \For{$\kk\in\left\{1,\ldots,\gLG\right\}$} {
   Computation of $\pikLG=\frac{\zLG_{\sommecolonne\kk}}{\n}$.
   }\tcp*[r]{\textcolor{mblue}{$\mathcal{O}(\n/\Sg)$}}
  \BlankLine
  \For{$\el\in\left\{1,\ldots,\mLG\right\}$} {
   Computation of $\rholLG=\frac{\wLG_{\sommecolonne\el}}{\dd}$.
   }\tcp*[r]{\textcolor{mblue}{$\mathcal{O}(\dd/\Sm)$}}
  \BlankLine
  Computation of $\balLG=\left(\bzLG\right)\transpose\bx\bwLG/\left[\pikLG\left(\rholLG\right)\transpose\right]\times\n\dd$. \tcp*[r]{\textcolor{mblue}{$\mathcal{O}\left(\n\dd\left[1/\Sg+1/\Sm\right]\right)$}}
  \BlankLine
  \BlankLine
  \KwOut{Numbers of clusters $\gLG$ and $\mLG$, matrices $\bzLG$ and $\bwLG$ and parameter $\bthetaLG$.}
  \BlankLine
}
  \caption{Algorithm $Largest$ $Gaps$. The complexity of each step is highlighted in blue.}
  \label{algo:LG}
\end{algorithm}

\paragraph{Estimator of $\bthetas$.} For the rest of the article, $\bthetas\Rev{=\left(\bpis,\brhos,\bals\right)}$ represents the parameters to be estimated (i.e. the parameters that were used to simulate the data). In the algorithm~\ref{algo:LG}, the estimator $\bthetaLG$ of $\bthetas$ is based on $\bzLG$ and $\bwLG$. $\pikLG$ (resp. $\rholLG$) is the proportion of class $\kk$ (resp. $\el$) in the partition $\zLG$ (resp. $\wLG$) and the estimator $\balLG$ is for all $\left(\kk,\el\right)\in\left\{1,\ldots,\gLG\right\}\times\left\{1,\ldots,\mLG\right\}$:
\[\alklLG=\frac{\sum_{\ii=1}^{\n}\sum_{\jj=1}^{\dd}\zikLG\wjlLG\xij}{\zsumkLG\wsumlLG}.\]

\paragraph{Complexity of the algorithm.}
On algorithm~1, the complexities of each step are added in blue at the end of the line. In the end, the \textit{LG} algorithm has a complexity of $\mathcal{O}\left(\max\left(\n\dd\left[\gLG+\mLG\right],\n\log\n,\dd\log\dd\right)\right)$.
As will be seen in the section~\ref{sec:consistency}, $\log\n$ is required to be much smaller than $\dd$ and $\log\dd$ much smaller than $\n$. In this case, the complexity is $\mathcal{O}\left(\n\dd\left[\gLG+\mLG\right]\right)$.\\
Moreover, $\sum_{\ii=2}^{\n}\Gi=\overline{\X_{\ligne(\n)}}-\overline{\X_{\ligne(1)}}$ being smaller than $1$ and for all $\kk\in\{1,\ldots,\gLG-1\}$, $G_{(\ii_\kk)}$ being greater than $\Sg$ then, in the worst case, $\gLG$ is smaller than $1/\Sg+1$.\\
As a conclusion, the complexity is $\mathcal{O}\left(\n\dd\left[1/\Sg+1/\Sm\right]\right)$ and, if the classification only is processed, the complexity is $\mathcal{O}\left(\n\dd\right)$.

\section{Identifiability of the model: assumptions on the model parameters}\label{sec:assumptions:model}

In this section, the conditions for the consistency of the estimators are explained. \Rev{For the rest of the article, $\gs$ and $\ms$ represent the true number of clusters in rows and columns and, as mentioned in the section~\ref{sec:algorithm}, $\bthetas=\left(\bpis,\brhos,\bals\right)$ corresponds to the true parameters . Moreover, $\bzs$ and $\bws$ are the unobserved partition matrices, resulting from the laws~$\mathcal{M}\left(1;\bpis\right)$ and~$\mathcal{M}\left(1;\brhos\right)$, used to simulate the data.}

\begin{nota}\textit{Key model parameters:}\\
Let us define $\pimin$ and $\rhomin$ the smallest class proportions:
\[\pimin=\underset{1\leq\kk\leq\gs}{\min}\!\piks\text{ and }\rhomin=\underset{1\leq\el\leq\ms}{\min}\!\rhols\]
and the smallest distance between any two conditional expectations of the normalized degrees (called model smallest gaps):
\begin{eqnarray*}
\deltapi=\underset{1\leq\kk\neq\kp\leq\gs}{\min}\;\left|\tauks-\taukps\right|
&\text{ and }&\deltarho=\underset{1\leq\el\neq\lp\leq\ms}{\min}\;\left|\xils-\xilps\right|\\
\end{eqnarray*}
where $\btaus=\Rev{\brhos\bals\transpose}$ and $\bxis=\Rev{\bpis \bals}$ are the proportions of the binomial distributions defined in Equation~\eqref{Eq:Concept}.
\end{nota}

Cluster identifiability by the \textit{LG} algorithm is possible under the following sufficient conditions:

\begin{assumpt} \textbf{(I)} The proportion of each row class (respectively column class) is non-negative, and conditional expected degrees of row (resp. column) clusters $(\tauk)_{1\leq \kk\leq \gs}$ (resp. $(\rhol)_{1\leq \el\leq \ms}$) are all distinct, which respectively amounts to:
\[ \pimin\rhomin>0 \text{ and } \deltapi\deltarho>0. \label{hyp:identifiability}\tag{I} \]
\end{assumpt}

These assumptions will be always made in the sequel of this article. They are equivalent to the sufficient conditions of identifiability of LBM given in \cite{keribin2014estimation} \Rev{but stronger than the conditions of \cite{brault2020consistency} or the identifiability conditions in the case where $(g,m)=(2,2)$ in the \cite{keribin2014estimation}.}

The first assumption, $\pimin\rhomin>0$, is classical in mixture models (for example, see \cite{brault2020consistency,iannario2010identifiability,maugis2011non}); without this assumption, the probability of having the right number of clusters is null. Thus the mixture model would be degenerated, as clusters with proportion zero would be always empty and the number of actually present clusters would be improper.

The second one, $\deltapi\deltarho>0$, ensures that not only the distribution of $\bx$ is identifiable, but also that of $(\Xiplus)_{1\leq i\leq \n}$, which is a marginal distribution of $\bx$. This is critical in order to recover all clusters with the \textit{LG} algorithm, which actually infers the clusters only from the set of variables $\{\Tid;\ii=1,\dots,\n\}$.
More precisely, it ensures that the distribution of $\Xiplus$ is a proper mixture with the same number of clusters as the distribution of $\X_{ij}$. For example, Assumption \eqref{hyp:identifiability} excludes the following typical case \Rev{(yet identifiable according to~\cite{keribin2014estimation})}:
\begin{eqnarray*}
 &&\g=2,\,m=2,\ \bpi=(1/2,1/2)\\
  &\text{ and }&\bal=\begin{pmatrix}
            a & b \\
            b & a
           \end{pmatrix}\text{ where } a,b\in]0,1[.\\
\end{eqnarray*}
Indeed, it gives $\tauu_1=\tauu_2=\frac{a+b}{2}$, and hence $\deltapi=0$. Conditionally on either $z_{i1}=1$ or $z_{i2}=1$, the distribution of $\Xiplus$ is the same: $\Bin\left(d,\frac{a+b}{2}\right)$, and therefore the two clusters cannot be distinguished just with the vector $(\Xiplus)_{1\leq i\leq \n}$. The distribution of $\Xiplus$ is not a proper mixture with two distinct clusters, it is a simple binomial distribution. Thus with our approach this model would be confused with the model $\g=1,\ \bpi=1,\ \bal=\frac{a+b}{2}$. Note that the set of parameters such that $\deltapi\deltarho=0$ \Rev{has zero Lebesgue measure}.

\section{Consistency}\label{sec:consistency}

This section presents the main result (Theorem \ref{th:consistency}), namely the consistency of estimators. Here, consistency means that the numbers of clusters and classifications are correct and that the distance between the estimates and the model parameters is smaller at any $t>0$ with a probability tending towards one, when the size of the matrix $(n,d)$ tends towards infinity. Before stating this theorem, we introduce some notations, especially related to the label switching problem.

\subsection{Distance on the parameters and the label switching issue}\label{Sec:distance}


For any two parameters $\btheta=(\bpi,\brho,\bal)$ with $(\g,\m)$ clusters and $\btheta'=(\bpi',\brho',\bal')$, with $(\g',\m')$ clusters, we define their distance as follows:
\begin{eqnarray*}
\Rev{\dHinf\left(\btheta,\btheta'\right)}
&=&\left\{\begin{array}{l}
\Rev{\norminf{\btheta-\btheta'}\text{ if }\g=\g',\ \m=\m'}\\
+\infty \text{ otherwise,}\\
\end{array}\right.\\
\end{eqnarray*}
where $\left\|\cdot\right\|_{\infty}$ denotes the norm defined for any $\by\in\mathbb{R}^{\g}$ by $\left\|\by\right\|_{\infty}=\max_{1\leq\kk\leq\g}\;\left|y_\kk\right|$.

We assume that two matrices $\bz,\bz'\in\mathcal{M}_{\n\times\g}\left(\{0,1\}\right)$ are equivalent, denoted $\bz\egalZ\bz'$, if there exists a permutation $\perms\in\mathfrak{S}\left(\left\{1,\ldots,\g\right\}\right)$ such that for all $\left(\ii,\kk\right)\in\{1,\ldots,\n\}\times\{1,\ldots,\g\}$, $\z_{\ii,\perms(\kk)}'=\zik$. By convention, we assume that two matrices with different numbers of columns are not equivalent. We introduce the similar notation $\egalW$ for the matrix $\bw$.\\
For all parameter $\btheta=(\bpi,\brho,\bal)$ with $(\g,\m)$ clusters and for all permutions \linebreak$\left(\perms,\permt\right)\in\mathfrak{S}\left(\left\{1,\ldots,\g\right\}\right)\times$ $\mathfrak{S}\left(\left\{1,\ldots,\m\right\}\right)$, we denote $\btheta^{\perms,\permt}=(\bpi^{\perms},\brho^{\permt},\bal^{\perms,\permt})$, by:
\begin{eqnarray*}
&\bpi^{\perms}=\left(\pii_{\perms(1)},\ldots,\pii_{\perms(\g)}\right),\;\;\brho^{\permt}=\left(\rhoo_{\permt(1)},\ldots,\rhoo_{\permt(\m)}\right)&\\
&\text{and }\bal^{\perms,\permt}=\left(\al_{\perms(1),\permt(1)},\al_{\perms(1),\permt(2)}, \ldots,\al_{\perms(1),\permt(\m)},\al_{\perms(2),\permt(1)},\ldots,\al_{\perms(\g),\permt(\m)}\right).&\\
\end{eqnarray*}
Like all mixture models, the LBM is affected by the label switching problem: clusters are \Rev{defined up to a permutation}. The algorithm can therefore find the right clusters but at a permutation of the labels. This also affects the parameters of the model, because the order of their coordinates depends on the labeling. The comparison of two classifications, and two parameter estimates, must then be done carefully: the distance between two parameter estimates must be calculated after permutation of their coordinates, using the permutation transforming the label allocation of the classification algorithm into the original label allocation of the model. Moreover, such a permutation exists and is unique when the classification is correct, i.e. when $\bzLG \egalZ \bzs$ (respectively $\bwLG \egalW \bws$). This permutation will thus be noted $\perms_{\clZ}$ (resp. $\permt_{\clW}$) on the event $\bzLG \egalZ \bzs\}$ (resp. $\bwLG \egalW \bws$). Thus, the consistency of the parameter estimators amounts to proving that the following quantity disappears in probability when $(n,d)$ tends to infinity:
\begin{eqnarray*}
\forall t>0,\Prob\left[\left.\dHinf\left(\bthetaLG^{\perms_{\clZ},\permt_{\clW}},\bthetas\right)>t\right|\bzLG \egalZ \bzs,\bwLG \egalW \bws\right]&
\underset{n,d\to+\infty}{\longrightarrow}&0. \\
\end{eqnarray*}
Outside of the event $\{\bzLG \egalZ \bzs\}$ (resp. $\bwLG \egalW \bws$), $\perms_{\clZ}$ (resp. $\permt_{\clW}$) will be defined as any arbitrary permutation in $\mathfrak{S}\left(\left\{1,\ldots,\gLG\right\}\right)$ (resp. $\mathfrak{S}\left(\left\{1,\ldots,\mLG\right\}\right)$), the identity for instance.
%

\subsection{A non-asymptotic upper bound}

This paragraph presents the main technical tool for obtaining consistency: a non-asymptotic upper bound from which the consistency results will be derived. 

\begin{theorem}[\Rev{Concentration inequality}]\label{thm:concentration}
Under identifiability assumptions \eqref{hyp:identifiability}, and if $\Sg\in ]0,\deltapi[ \text{ and } \Sm \in ]0,\deltarho[$ for $\n,\dd$ large enough, then we have for all $t>0$:
\begin{eqnarray*}
&&\!\!\!\!\!\!\!\!\Prob\left(\gLG\neq\gs\text{ or }\mLG\neq\ms\text{ or }\bzLG\cegalZ\bzs\text{ or }\bwLG\cegalW\bws\text{ or }\dHinf\left(\permbthetaLG,\bthetas\right)>t\right)\\
&\leq& 2\n \exp\left( -\frac{\dd}{2} \min(\deltapi-\Sg,\Sg)^2 \right) +\gs\left(1-\pimin\right)^\n\\
& &+ 2\dd \exp\left( -\frac{\n}{2} \min(\deltarho-\Sm,\Sm)^2 \right) +\ms\left(1-\rhomin\right)^\dd \\
& &+\Rev{2\gs\ms\left[1-\pimin\rhomin\left(1-e^{-2t^2}\right)\right]^{nd}}+ 2\gs e^{-2\n t^2}+2\ms e^{-2\dd t^2}.\\
\end{eqnarray*}
\end{theorem}

The proof (in Appendix \ref{Appendix:thm:concentration}) is made in two steps, emphasizing the originality of the method in comparison with EM-like algorithms: here the classification is completely done first, and parameters are then estimated afterwards. Thus an upper bound on classifications and selection of class numbers will be first established (Proposition \ref{prop:concentration_inequality:g}), and secondly an upper bound on the parameter estimators, given that both classifications and class numbers are right (Proposition \ref{prop:btheta}).

\Rev{If $\n$ and $\dd$ increase at the same speed, the larger the size, the smaller the terminal. We find the importance that $\Sg$ and $\Sm$ are small enough to detect the real jumps but not so small that there are only them. Moreover, the larger $\pimin$ and $\rhomin$ are, the more likely it is to have at least one representative of each cluster.}

\Rev{In addition to showing the consistency of the \textit{LG} algorithm's estimators (see section~\ref{sec:consistency}), the bound of theorem~\ref{thm:concentration} quantifies how (very) easy a configuration is: given a configuration $\left(\btheta,\n,\dd\right)$, it is possible to estimate the probability that the algorithm \textit{Largest Gaps} finds the right configuration. In the case of a configuration obtained by another algorithm, this bound can be calculated and if it is (very) small, it means that a study of the marginals would have been sufficient; if the \textit{Largest Gaps} algorithm does not find the same results at all, a doubt can be cast on the results (estimation or relevance of the use of a LBM) by precaution.} 

\subsection{Consistency of the estimators}\label{sec:consistency}

The following theorem provides sufficient assumptions on the threshold sequences $(\Sgnd,\Smnd)_{n,d}$ to ensure the consistency of the inference method based on the \textit{LG} algorithm when $\n,\dd$ tend to infinity. Note that this result is therefore asymptotic only, and it does not provide any guarantee for fixed $(\n,\dd)$. Nevertheless the rates can be used as a suggestion to choose the thresholds, as  even though it must be carefully used, because it is anyway impossible to know whether the asymptotic regime is reached or not.

\begin{theorem}[\Rev{Consistency}] \label{th:consistency}
Under identifiability assumptions \eqref{hyp:identifiability} and the following assumptions:
 \[
\begin{array}{rl}
\Sgnd\underset{\n,\dd\to+\infty}{\longrightarrow}0,& \ \Smnd\underset{\n,\dd\to+\infty}{\longrightarrow}0, \\
\underset{\n,\dd\to+\infty}{\underline{\lim}}\Sgnd\sqrt{\frac{\dd}{\log\n}}>\sqrt{2}&
\text{ and } \underset{\n,\dd\to+\infty}{\underline{\lim}}\Smnd\sqrt{\frac{\n}{\log\dd}}>\sqrt{2},\\
\end{array}\label{hyp:vanishing} \tag{TA} 
\]
where $\underline{\lim}$ is the lower limit; then classifications, model selection and estimators are consistent, that is, for all $t>0$:
 \begin{eqnarray*}
&&\Prob\left(\gLG\neq\gs\text{ or }\mLG\neq\ms\text{ or }\bzLG\cegalZ\bzs\text{ or }\bwLG\cegalW\bws\right.\\
&&\left.\text{ or }\dHinf\left(\bthetaLG^{\perms_{\clZ},\permt_{\clW}},\bthetas\right)>t\right)\underset{\n,\dd\to+\infty}{\longrightarrow}0.\\
\end{eqnarray*}
\end{theorem}

The proof of this result is available in Appendix~\ref{cor:Appendix:consistency}.

The $LG$ algorithm has two input parameters, the $(\Sg,\Sm)$ thresholds, which must be set correctly to discover all clusters. Recall that the purpose of gap thresholding is to distinguish between external and internal gaps (see comments in section~\ref{subsec:concept}). First, the thresholds must be smaller than the smallest gaps in the model, $\Sgnd<\deltapi$ and $\Smnd<\deltarho$, for $\n$ and $\dd$ sufficiently large; otherwise, some clusters will consist of mixed clusters. Since $\deltapi$ and $\deltarho$ are not known a priori, we assume that the thresholds decrease with respect to $\n$ and $\dd$, to ensure that they are asymptotically small enough. On the other hand, if the threshold sequences decrease too fast, the thresholds will be asymptotically too small and at least one class will be split into several clusters by the algorithm. More technically, the convergence rate given in the theorem guarantees that the upper bound of the theorem \ref{thm:concentration} tends to 0 which implies consistency. If the sequences disappear faster, coherence is no longer guaranteed.

\Rev{Consistency can thus be obtained whatever the sequences $\Sgnd$ and $\Smnd$ taken provided that they verify the conditions~\eqref{hyp:vanishing}; nevertheless, if we couple with the condition that $\Sgnd<\deltapi$ and $\Smnd<\deltarho$ of Theorem~\ref{thm:concentration}, we see that it is preferable that the thresholds decrease rapidly. For example, the threshold $\Sgnd=\sqrt{2\log(\n)/\dd}(1+\varepsilon)$ with $\varepsilon>0$ can be used. Moreover, once a first configuration has been found, it is possible to estimate $\deltapi$ and $\deltarho$ and to re-estimate the partitions with the thresholds $\widehat{\deltapi}/2$ and $\widehat{\deltarho}/2$ in order to see if a better configuration is obtained.}

\begin{remark}\label{rem:consistency}
The assumption \eqref{hyp:vanishing} of the theorem implies that $\n/\log\dd$ and $\dd/\log\n$ tend to $+\infty$~; this assumption is found in the results on the consistency of the maximum likelihood estimator (see~\cite{brault2020consistency}). Therefore, $\bx$ is allowed to have an oblong shape. For example, $\dd=\n^{\gamma}$ with $\gamma>0$ satisfies the assumption.
\end{remark}

\begin{remark}\label{rem:symmetric}
The roles of the rows and columns are symmetric: the result remains true if the transpose of $\bX$ is studied rather than the matrix $\bX$.
\end{remark}

\subsection{Consistency when model parameters vary}

Finally varying model parameters are also considered in this paragraph. Indeed when $\n$ and $\dd$ increase, it can be reasonable to assume that new clusters arise: in this paragraph $\gsnd$ and $\msnd$ are hence assumed to be growing to infinity when $\n,\dd$ tend to $+\infty$. The consequence of this assumption is the convergence to zero of both the proportions of the smallest clusters and the model smallest gaps: $\pimin,\rhomin,\deltapi,\deltarho$ tend to 0 when 
$\n,\dd$ tend to infinity. \Rev{For example, since the parameters $\btaus=\left(\tauks\right)_{1\leq k\leq \gs}$ are probabilities, the following inequalities are obtained:}
\Rev{\[ (\gs-1)\deltapi \leq \sum_{k=1}^{\gs -1} (\tau_{(k+1)}^{\vrai} - \tau_{(k)}^{\vrai})= \tau_{(\gs)}^{\vrai}-\tau_{(1)}^{\vrai} \leq 1\]
\[\text{and }(\ms-1)\deltarho \leq \sum_{l=1}^{\ms -1} (\xii_{(l+1)}^{\vrai}-\xii_{(l)}^{\vrai}) \leq 1, \]}

\noindent \Rev{where $\tau_{(1)}^{\vrai}<\dots<\tau_{(\gs)}^{\vrai}$ (resp. $\xii_{(1)}^{\vrai}<\dots<\xii_{(\ms)}^{\vrai}$) are the $(\tauks)_{1\leq k\leq \gs}$ (resp. $(\xils)_{1\leq l\leq \ms}$)} sorted in increasing order. \Rev{In the other side, $\gs$ is bounded by $1/\pimin$ by the following inequality:
\[\gs\pimin\leq\sum_{k=1}^{\gs}\piks=1.\]}
In particular, it is possible that one of these parameters tend to zero while the number of groups does not change (for example, in the case of increasingly sparse matrices). In this framework, admissible convergence rates of the model parameters are provided, as well as the corresponding admissible convergence rates of the thresholds. It thus tells how robust the consistency is.

\begin{theorem}[\Rev{Consistency in sparse case}]\label{th:consistency:spare}
 Under the assumptions of the previous theorem (\eqref{hyp:identifiability} and \eqref{hyp:vanishing}) and the following additional assumptions :
 \begin{itemize}
 	\item Assumptions on $\deltapind$ and $\Sgnd$ (resp. $\deltarhond$ and $\Smnd$):
 	\[\underset{\n,\dd\to+\infty}{\underline{\lim}}\frac{\deltapind}{\Sgnd}>2,  \text{ and }\underset{\n,\dd\to+\infty}{\underline{\lim}}\frac{\deltarhond}{\Smnd}>2. \label{hyp:varying1} \tag{MA.1}\]
 	\item Assumptions on $\piminnd$ and $\rhominnd$:
 	\[\Rev{\n\piminnd\underset{\n,\dd\to+\infty}{\longrightarrow}+\infty\quad\text{ and }\quad\dd\rhominnd\underset{\n,\dd\to+\infty}{\longrightarrow}+\infty} \label{hyp:varying2} \tag{MA.2}\]
 \end{itemize}
then classifications, model selection and estimators are also consistent.
\end{theorem}

The proof of this result is available in Appendix~\ref{cor:Appendix:consistency:spare}.

\Rev{Theorem~\ref{th:consistency:spare} shows that the estimates remain consistent even with sparse matrices. On the other hand, for an optimal result and if a number of non-zero cells is fixed, the estimation will be all the easier if the latter are well distributed between the classes, for example, in a staircase shape.}

\begin{remark}
    \Rev{As $\piminnd\leq 1/\gsnd$ and  $\rhominnd\leq 1/\msnd$, the assumptions~\eqref{hyp:varying2} implie that
    \[\frac{\n\dd}{N_{\text{block}}^{\star\,\n\dd}}\underset{\n,\dd\to+\infty}{\longrightarrow}+\infty\]
    where $N_{\text{block}}^{\vrai\,\n\dd}=\gsnd\msnd$ is the number of blocks. In particular, the number of blocks can therefore increase with $\n$ and $\dd$ but not too quickly.}
\end{remark}

\section{Simulations}\label{sec:simulation}

In this section, simulations to test the quality of the bounds obtained in the previous theorems \Rev{and to compare the computing times with classical procedures are presented}\footnote{\Rev{For the sake of reproducibility, all codes are available on the following Gitlab: \url{https://gricad-gitlab.univ-grenoble-alpes.fr/braultv/largest-gaps}}}.

\subsection{Estimation of the number of clusters}\label{sec:simu:g}

We use an experimental design to illustrate the results of Theorem~\ref{th:consistency}. As the number of row (resp column) clusters is the basis of the other estimations, this is the only parameter studied in this section. The experimental design is defined with $\gs=5$ and $\ms=4$ and the following parameters
\[\bals=\begin{pmatrix}
    \varepsilon&\varepsilon&\varepsilon&\varepsilon\\
    1-\varepsilon&\varepsilon&\varepsilon&\varepsilon\\
    1-\varepsilon&1-\varepsilon&\varepsilon&\varepsilon\\
    1-\varepsilon&1-\varepsilon&1-\varepsilon&\varepsilon\\
    1-\varepsilon&1-\varepsilon&1-\varepsilon&1-\varepsilon\\
    \end{pmatrix}\]
with $\varepsilon\in\{0.05,0.1,0.15,0.2,0.25,0.3\}$. For $\varepsilon=0.05$, the values of the Bernoulli parameters are $0.05$ (resp. $0.95$) and the associated blocks contain essentially white (resp. black) cases. At the opposite, for $\varepsilon=0.25$, the block are more homogeneous and more difficult to distinguish. For the class proportions, we suppose two possibilities
\begin{itemize}
\item Balanced proportions:
\[\Rev{\bpis=\begin{pmatrix}0.2&0.2&0.2&0.2&0.2\\\end{pmatrix}\;\text{ and }\;
\brhos=\begin{pmatrix}0.25&0.25&0.25&0.25\\\end{pmatrix}}\]
with the following parameters
\[\pimin=0.2\text{ and }\deltapi=0.25-0.5\varepsilon.\]
\item Arithmetic proportions:
\[\Rev{\bpis=\begin{pmatrix}0.1&0.15&0.2&0.25&0.3\\\end{pmatrix}\;\text{ and }\;
\brhos=\begin{pmatrix}0.1&0.2&0.3&0.4\\\end{pmatrix}}\]
with the following parameters
\[\pimin=0.1\text{ and }\deltapi=0.1-0.2\varepsilon.\]
\end{itemize}
The number of rows $\n$ and the number of columns $\dd$ vary between 20 and 4000 by step 20 and for each configuration, \Rev{100} matrices were simulated. For the choice of the thresholds $\Sg$, we studied four cases:

\begin{enumerate}
 \item We first propose a constant oracle threshold to illustrate the Theorem~\ref{thm:concentration}, which suggests that any constant threshold strictly between $0$ and $\deltapi$ gives the consistency:
 \[\Sun=\deltapi/2.\]
 \item In practice, there is mostly no reason why the parameter $\deltapi$ could be known. In this case, Theorem~\ref{th:consistency} claims that we need to use a varying threshold $\Sgnd$ instead, such that $\Sgnd$ tends to 0 but not too fast (slower than $\sqrt{2\log\n/\dd}$). If the threshold decreases too slowly, it may be larger than $\deltapi$ and the smallest gaps could be undetected. On the opposite, if the threshold decreases too fast, we may detect too many gaps. In the simulation, we studied three possibilities:
 \begin{enumerate}
    \item Faster threshold:\newline $\Sdeux=\sqrt{2\log(\n)/\dd}\left(1+10^{-10}\right)$.
    \item Middle threshold: $\Strois=2\sqrt{2\log(\n)/\dd}$.
    \item Slower threshold: $\Squatre=\left(\log(\n)/\dd\right)^{1/4}$.
\end{enumerate}
\end{enumerate} 

\begin{figure}[!h]
\begin{center}
    \includegraphics[width=\linewidth]{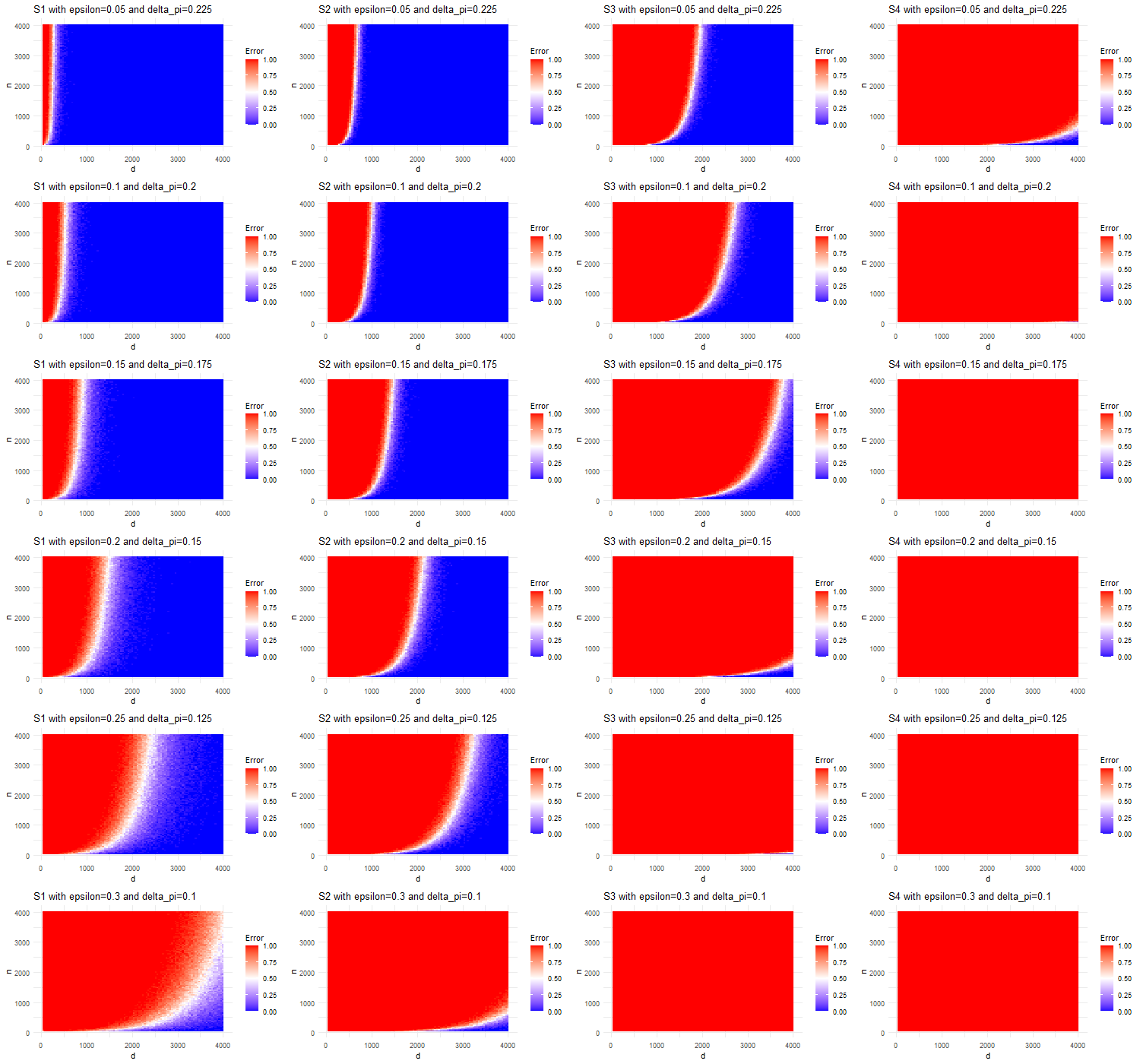}
\caption{\label{Fig:Res:Prop}Proportions of bad estimations of $\gs$ for the parameter $\varepsilon\in\{0.05,\ldots,0.3\}$ (in rows) following the thresholds (in columns) used for the balanced case: for each graphic, the number of rows $\n$ \Rev{(y-axis)} and the number of columns $\dd$ \Rev{(x-axis)} varies between 20 and 4000.}
\end{center}
\end{figure}

\begin{figure}[!h]
\begin{center}
    \includegraphics[width=\linewidth]{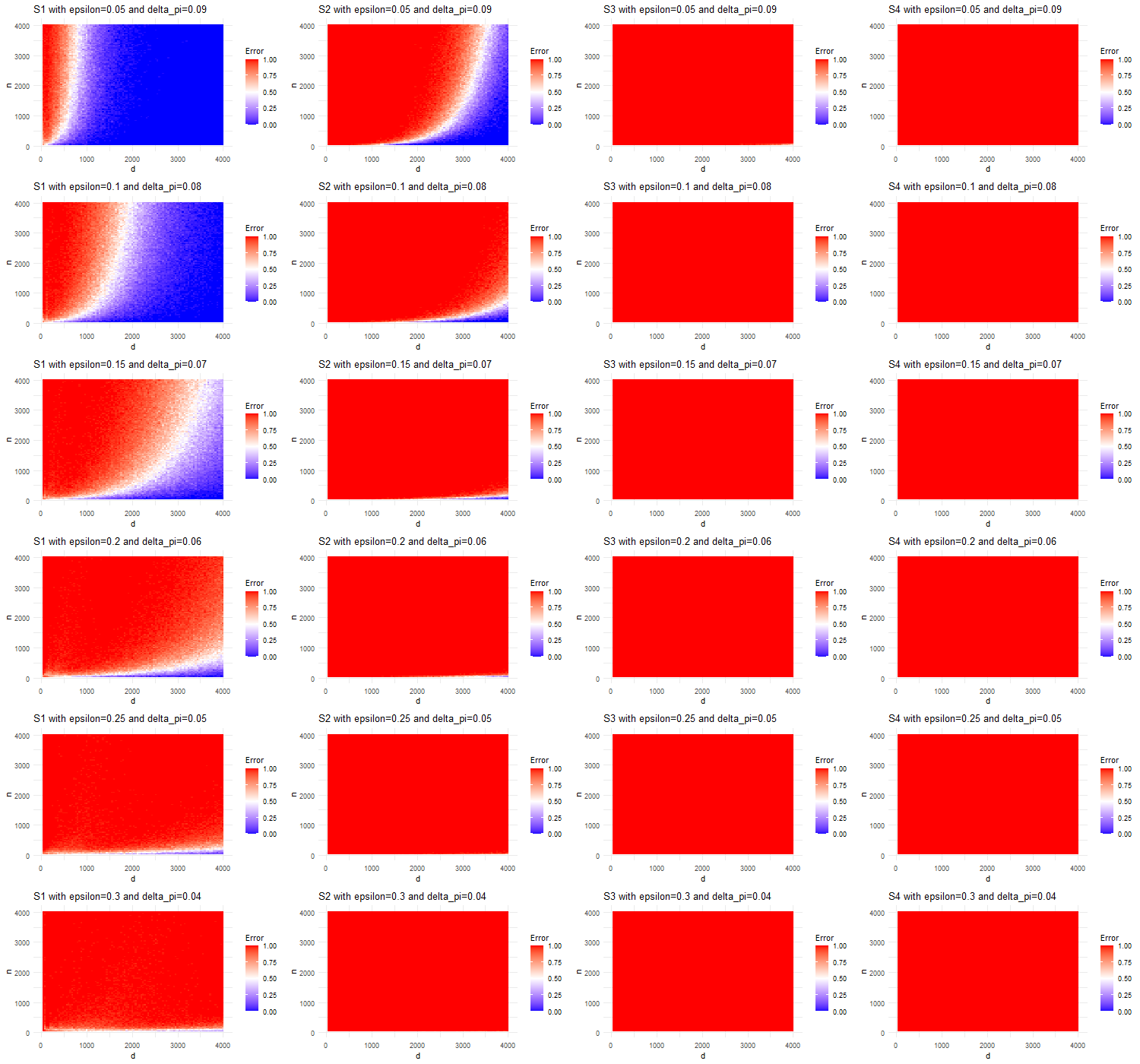}
\caption{\label{Fig:Res:NProp}Proportions of bad estimations of $\gs$ for the parameter $\varepsilon\in\{0.05,\ldots,0.3\}$ (in rows) following the thresholds (in columns) used for the arithmetic case: for each graphic, the number of rows $\n$ \Rev{(y-axis)} and the number of columns $\dd$ \Rev{(x-axis)} varies between 20 and 4000.}
\end{center}
\end{figure}

Figures~\ref{Fig:Res:Prop} and~\ref{Fig:Res:NProp} show the proportions of bad estimations of $\gs$ \Rev{following the parameters $\varepsilon$ (in rows)} as function of the number of rows $\n$, the numbers of columns $\dd$ and the thresholds used \Rev{(columns)}. For each figure, the number of columns ($\dd$) increases \Rev{on the x-axis} and the number of rows ($\n$) increases \Rev{on the x-axis}. The red color corresponds to a bad estimate while the blue color corresponds to a good estimate (the more blue it is, the better).

As expected, it appears that the best threshold is the oracle $\Sun=\deltapi/2$ but this threshold cannot be used in practice because $\deltapi$ is unknown. For the scaled thresholds, $\Sdeux=\sqrt{2\log\n/\dd}\left(1+10^{-10}\right)$ is the best.

We can see that the larger the number of rows $\n$ is, the worse the estimation is and the larger the number of columns $\dd$ is, the better the estimation is. In the case of $\n=\dd$, the quality of the estimation increases with $\n$. In particular, the model selection can be generalized for the case of \cite{channarond2011Classification} and the results would be similar. $\pimin$ has a weak effect because it is rare to have an empty class but the effect of $\deltapi$ seems to be greater.

\subsection{Estimation of the parameters}

\begin{figure}[!h]
\begin{center}
    \includegraphics[width=\linewidth]{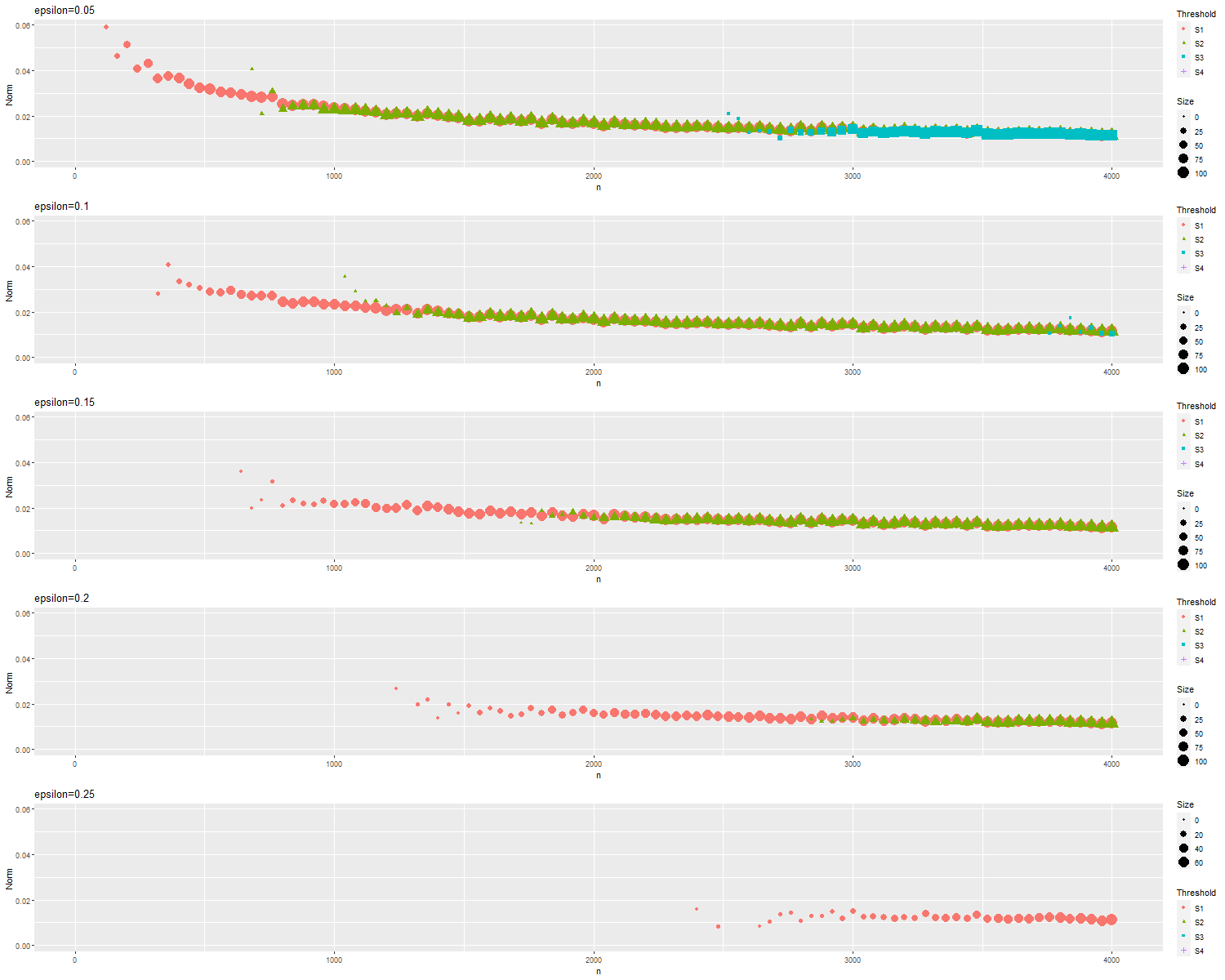}
\caption{\label{Fig:Res:Theta}Average estimate of the distance $\dHinf\left(\bthetas,\bthetaLG\right)$ between the true parameters and the estimated parameters following $\varepsilon$ (rows) in the balanced case: for each graphic, the threshold is represented by different symbols (\textcolor{RubineRed}{$\bullet$} for $\Sun$, \textcolor{LimeGreen}{$\blacktriangle$} for $\Sdeux$, \textcolor{Cerulean}{$\blacksquare$} for $\Strois$ and  \textcolor{purple}{$\boldsymbol{+}$} for $\Strois$) and the size for the number of finite values used; the number of rows $\n$ varies between 40 and 4000 and $\dd$ is supposed equal $\n$.}
\end{center}
\end{figure}

To illustrate the convergence of the the estimation $\bthetaLG$ to the true parameter $\bthetas$, the same experimental design is chosen with the number of rows $\n$ and the number of columns $\dd$ vary between 40 and 4000 by step 40 and, for each configuration, 100 matrices were simulated. To estimate the quality of the estimation, the distance $\dHinf$ introduced in Section~\ref{Sec:distance} is calculated. As the distance equals $+\infty$ when the number of clusters of $\bthetaLG$ is different of its of $\bthetas$, we chose to represent the mean of the finite values with the size of the point corresponding at the number of finite values (see Figure~\ref{Fig:Res:Theta} for \Rev{the balanced case} and Figure~\ref{Fig:Res:Theta:Neq} \Rev{for the arithmetic case}).

We observe that the error decreases with the number of observations and the results are identical when the numbers of clusters are correctly estimated. The optimization of the \textit{LG} algorithm thus depends on the choice of the threshold. In particular, we find that the oracle threshold ($S_1$) finds the right number of clusters faster than the evolutionary thresholds. \Rev{In particular,the threshold $\Squatre$ does not appear on the graphic because it underestimates the number of blocks}.

\subsection{\Rev{Comparison of computing times}}

\Rev{To conclude the part of the simulations, the computation times are estimated and compared with a classical procedure. For that, the plan presented previously is taken again by taking only the square matrices (n=d) with $n\in\{100,200,\ldots,900\}$ for the balanced case and $n\in\{100,200,\ldots,700\}$ for the arithmetic. Six procedures are studied:
\begin{itemize}
    \item the algorithm $LG$ with the four previous threshold: $S_1$, $\Sdeux$, $\Strois$ and $\Squatre$;
    \item the algorithm \textit{variational Bayes} with the hyperparameters $(4,4)$ as proposed by~\cite{keribin2014estimation} knowing the true number of parameters $\gs$ and $\ms$ (named \textit{Simple VBayes} for the next);
    \item as the number of blocks is usually unknown, the combination of the algorithm \textit{variational Bayes} and the \textit{Integrated Complete-data Likelihood} (see \cite{biernacki2000assessing,keribin2014estimation}) is studied with the hyperparameters $(4,4)$ on a $2\times7$ square grid (named \textit{VBayes+ICL} for the next).
\end{itemize}
For the implementation, the \texttt{R} package \texttt{blockcluster} (version 4.5.1; see~\cite{parmeet2017blockcluster}) with the function \texttt{coclusterBinary} is used. For each configuration, 20 matrices are simulated and procedures are evaluated on two criteria:}
\begin{itemize}
    \item \Rev{the mean computation times with 10 runs; for this, the package \texttt{microbenchmark} (version 1.4.9; see~\cite{olaf2021microbenchmark}) with the function \texttt{microbenchmark} is used and plans are launched on a cluster}\footnote{\Rev{Cluster \texttt{Luke44}, 28 cores, 128Go RAM, GPU 2xK40m, 2xIntel Xeon E5-2680 2.40 GHz; more informations on the url} \url{https://scalde.gricad-pages.univ-grenoble-alpes.fr/web/pages/presentation.html}}.
    \item \Rev{the estimation quality of the number of clusters per row after an estimation of each matrix (except for \textit{Simple VBayes} as the number is assumed be known); results are averaged over the 20 matrices.}
\end{itemize}

\begin{figure}[!h]
\begin{center}
    \includegraphics[width=\linewidth]{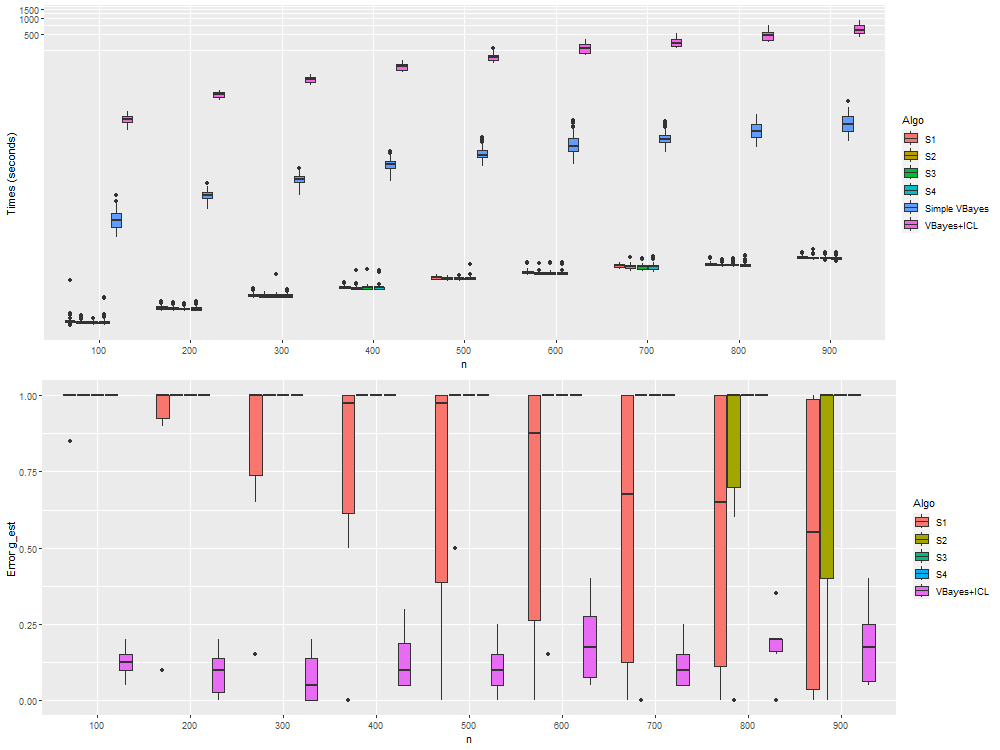}
\caption{\label{Fig:timesGlobal:Eq}\Rev{On the top, boxplots of the computing times (in seconds; logarithmic scale) for each procedure (colour) in function of the numbers of rows and columns ($\n=\dd$) over the 1200 simulations (200 per $\varepsilon$) for the balanced case. On the bottom, boxplots of quality of estimations averaged over the 20 matrices for each $\varepsilon$ for each procedure (colour) in function of the numbers of rows and columns ($\n=\dd$) for the balanced case.}}
\end{center}
\end{figure}

\begin{figure}[!h]
\begin{center}
    \includegraphics[width=\linewidth]{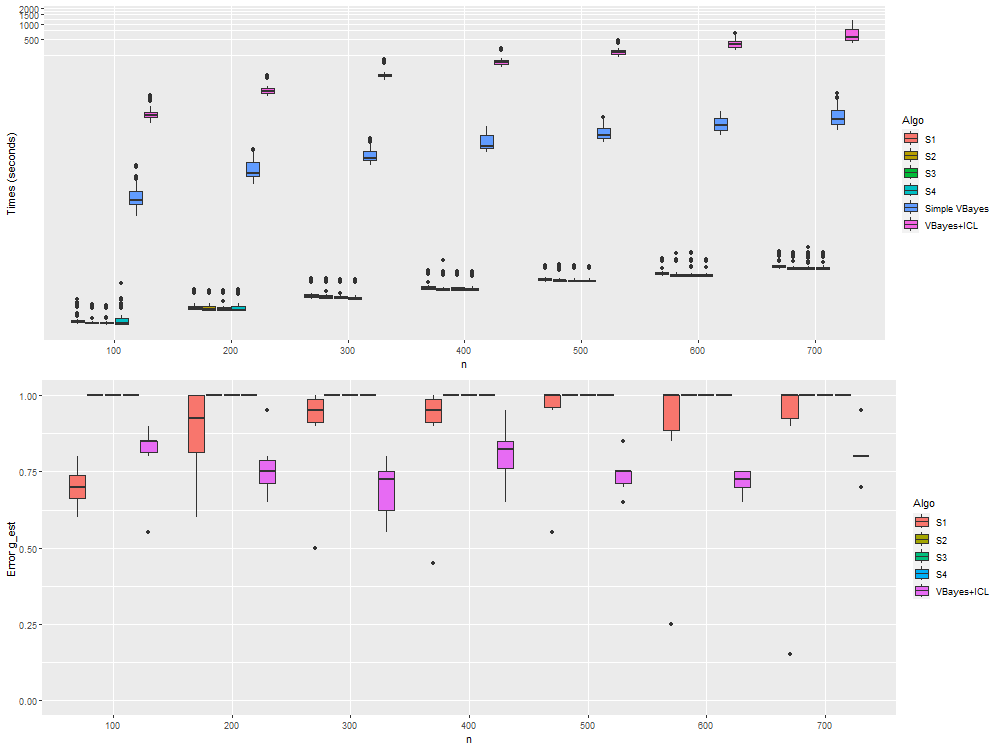}
\caption{\label{Fig:timesGlobal:Neq}\Rev{On the top, boxplots of the computing times (in seconds; logarithmic scale) for each procedure (colour) in function of the numbers of rows and columns ($\n=\dd$) over the 1200 simulations (200 per $\varepsilon$) for the arithmetic case. On the bottom, boxplots of quality of estimations averaged over the 20 matrices for each $\varepsilon$ for each procedure (colour) in function of the numbers of rows and columns ($\n=\dd$) for the balanced case.}}
\end{center}
\end{figure}

\Rev{The results are displayed in the Figures~\ref{Fig:timesGlobal:Eq} (balanced case) and~\ref{Fig:timesGlobal:Neq} (arithmetic case) where the computation times (in seconds with logarithmic scale; on the top) and the quality of the estimations (averaged for each $\varepsilon$ on over the 20 matrices; on the bottom) have been grouped by matrix size (regardless of the values of $\varepsilon$). A detailed version for each $\varepsilon$ is available on the Figures~\ref{Fig:times:Eq} and~\ref{Fig:times:Neq} in supplementary material. Moreover, the means and standard deviations (in milliseconds) are represented in the Tables~\ref{tab:times:Eq} and~\ref{tab:times:Neq}.
}

\begin{table}[!h]
    \centering
    \caption{\Rev{Computing times for the different procedures in columns (\textit{LG} with the four threshold, \textit{VBayes} with the true number of clusters and \textit{VBayes} coupled with the \textit{ICL} criterion on a $2\times7$ square grid) following the number of rows and columns (in rows) for the balanced case. Each cell represents the average (in milliseconds) and the standard deviation (in parenthesis) over the 1200 simulations (200 per $\varepsilon$).}}
    \label{tab:times:Eq}
\begin{tabular}{|c|*{6}{|c}|}
     \cline{2-7}
     \multicolumn{1}{c|}{}&$S_1$&$S_2^{\n,\dd}$&$S_3^{\n,\dd}$&$S_4^{\n,\dd}$&Simple VBayes&VBayes+ICL\\\cline{2-7}\hline
\multirow{2}{*}{n=100}& 2&2&2&2&171&12765\\
&(0.9)&(0.1)&(0.1)&(0.8)&(71.5)&(2180.1)\\\hline
\multirow{2}{*}{n=200}& 3&3&3&3&468&37193\\
&(0.4)&(0.3)&(0.3)&(0.3)&(102)&(4385.4)\\\hline
\multirow{2}{*}{n=300}& 6&6&6&6&930&71594\\
&(0.5)&(0.4)&(0.9)&(0.4)&(204.1)&(8711.8)\\\hline
\multirow{2}{*}{n=400}& 8&8&8&8&1789&123632\\
&(0.6)&(1)&(1)&(1.4)&(461.6)&(18941.2)\\\hline
\multirow{2}{*}{n=500}& 13&12&12&12&2934&187319\\
&(0.8)&(0.6)&(0.6)&(1.2)&(846.2)&(29995)\\\hline
\multirow{2}{*}{n=600}& 16&16&16&16&4597&280348\\
&(1.4)&(1)&(1.8)&(2.1)&(2044)&(57861.3)\\\hline
\multirow{2}{*}{n=700}& 21&21&21&21&5724&356932\\
&(1.5)&(2)&(2)&(2.6)&(1718.6)&(61458.8)\\\hline
\multirow{2}{*}{n=800}& 23&22&22&22&8050&490034\\
&(1.6)&(1.3)&(1.4)&(1.9)&(2900.3)&(107033.6)\\\hline
\multirow{2}{*}{n=900}& 31&30&30&30&11057&629238\\
&(1.9)&(2.3)&(1.9)&(1.8)&(4469.8)&(136921.4)\\\hline
\end{tabular}
\end{table}

\begin{table}[!h]
    \centering
    \caption{\Rev{Computing times for the different procedures in columns (\textit{LG} with the four threshold, \textit{VBayes} with the true number of clusters and \textit{VBayes} coupled with the \textit{ICL} criterion on a $2\times7$ square grid) following the number of rows and columns (in rows) for the arithmetic case. Each cell represents the average (in milliseconds) and the standard deviation (in parenthesis) over the 120 simulations (20 per $\varepsilon$).}}
    \label{tab:times:Neq}
\begin{tabular}{|c|*{6}{|c}|}
     \cline{2-7}
     \multicolumn{1}{c|}{}&$S_1$&$S_2^{\n,\dd}$&$S_3^{\n,\dd}$&$S_4^{\n,\dd}$&Simple VBayes&VBayes+ICL\\\cline{2-7}\hline
\multirow{2}{*}{n=100}& 2&2&2&2&530&20819\\
&(0.9)&(0.7)&(0.7)&(1.2)&(302.3)&(7768.3)\\\hline
\multirow{2}{*}{n=200}& 4&4&4&4&1677&58494\\
&(1.4)&(1.4)&(1.4)&(1.4)&(693.7)&(18294.9)\\\hline
\multirow{2}{*}{n=300}& 7&7&7&7&3183&117450\\
&(2.4)&(2.4)&(2.4)&(2.4)&(1145.8)&(34865.2)\\\hline
\multirow{2}{*}{n=400}& 10&10&10&9&5637&203479\\
&(3.2)&(3.7)&(3.2)&(3.3)&(2119)&(55356.2)\\\hline
\multirow{2}{*}{n=500}& 14&14&13&13&8503&304668\\
&(3.9)&(4)&(3.9)&(3.9)&(2862.3)&(76486.8)\\\hline
\multirow{2}{*}{n=600}& 18&18&18&18&12423&435626\\
&(5.2)&(5.5)&(5.9)&(5.3)&(3889.1)&(100281.1)\\\hline
\multirow{2}{*}{n=700}& 25&24&24&24&17794&622631\\
&(6.5)&(6.6)&(7)&(6.8)&(7596.6)&(180216.9)
\\\hline
\end{tabular}
\end{table}

\Rev{On the quality, the procedure \textit{VBayes+ICL} estimates better than the other the number of clusters in rows and the differences between the four threshold are the same than the section~\ref{sec:simu:g}.}

\Rev{For the computing times, the order of \textit{Largest Gaps} is that of the millisecond while that of \textit{VBayes+ICL} is of the second (and the minute after $n=300$). Moreover, the computing time for the \textit{Largest Gaps} seems to increase linearly with $\n$ (the side length of the square matrix) as stated in the section~\ref{sec:algorithm} while \textit{VBayes+ICL} increases much faster than linearly. Moreover, \textit{Largest Gaps} does not seem to be impacted by the configuration while the convergence time of \textit{VBayes+ICL} (and even \textit{Simple VBayes} who knows the right number of clusters) is larger for the case of arithmetic proportions. And finally, the computation time of \textit{VBayes+ICL} depends on the maximum choice of the number of clusters (fixed here quite close to the true result) but unknown in practice while \textit{Largest Gaps} does not need this kind of parameter and the time is almost independent of the number of blocks.}

\section{Conclusion \Rev{and discussion}}\label{sec:conclusion}

First of all, \textit{Largest Gaps} is a co-clustering algorithm, which has nice theoretical properties: its computational cost is much lower than most known algorithms and it provides a consistent procedure under the Latent Block Model, for all inference tasks: model selection, classification and estimation of the model parameters. Since the algorithm is simple, the consistency is rather easy to obtain. Note that in this article, only binary matrices have been studied, but the model as well as the method and the proofs can be directly extended for distributions which have the same concentration properties, for example compactly-supported distributions, where the support is known.

As a consistent algorithm, the advantage of the \textit{LG} algorithm is the simplicity: it is a simple and original way to exploit the concentration of the marginal distributions of the matrix $\bx$ under LBM. But it lacks robustness: a single outlier in the marginal distribution can affect the whole classification, which is not desirable by the use on real-world datasets. Moreover, the procedure is essentially based on the choice of threshold and simulations have shown that there is still room for improvement in this choice. Many other algorithms, also based on this nice feature, could be used instead and might be better in practice, but harder to analyze from a theoretical point of view. For example, fitting a binomial mixture model on the variables $\Tund,\dots\Tnd$ with an EM-algorithm, could be fruitful as well.

The contribution of this article actually goes beyond the \textit{LG} algorithm. It shows special features of the latent structure of the LBM, and their consequences. In particular, when the asymptotic regime of the model is reached, the latent structure is almost obvious, and moreover it can be \Rev{pick either one but not the two} from a summary of the data. Indeed even basic algorithms like \textit{LG} can retrieve the latent clusters from variables $\Tund,\dots\Tnd$ (sums of the rows and of the columns of the matrix $\bx$), whereas most known classification algorithms are unusable, because of their complexity.

A consequence of this remark is that the LBM should be used sparingly on very large real data sets: if the number of clusters requested is too small and if the blocks are sufficiently contrasted then the marginals should highlight the clusters. If this is not the case, the use of the LBM should be questioned. The bound of the theorem~\ref{thm:concentration} would give an estimate of the quality of the use of this model.


Finally, it appears in the simulations that the estimate of the number of clusters is underestimated. Moreover, the distribution of marginals from real-world data rarely has such obvious deviations as assumed by the model asymptotic. To overcome this problem, it would be interesting to estimate the row clusters with a mixture model on the variables $(\Tund,\ldots,\Tnd)$; this will be the subject of future work.


\FloatBarrier

\appendix

\section{Proof of Theorem \ref{thm:concentration}} \label{Appendix:thm:concentration}

In this appendix, we present demonstrations of the concentration inequality. We first estimate the probability of having the right number of clusters and then the right classification knowing that we have the right number of clusters. Finally, we evaluate the quality of the parameter estimates.





First of all, note that $\{\bzLG \egalZ \bzs\} \subset \{\gLG = \gs\}$ and $\{\bwLG \egalW \bws\} \subset \{\mLG = \ms\}$, hence :

\begin{eqnarray*}
&&\Prob\left(\gLG\neq\gs\text{ or }\mLG\neq\ms\text{ or }\bzLG\cegalZ\bzs\text{ or }\bwLG\cegalW\bws\text{ or }\dHinf\left(\permbthetaLG,\bthetas\right)>t\right)\\
&=&\Prob\left(\bzLG\cegalZ\bzs\text{ or }\bwLG\cegalW\bws\text{ or }\dHinf\left(\permbthetaLG,\bthetas\right)>t\right)\\
&=& \Prob\left(\bzLG\cegalZ\bzs\text{ or }\bwLG\cegalW\bws \right)+ \Prob\left(\left\{\dHinf\left(\permbthetaLG,\bthetas\right)>t\right\} \left\backslash \left\{\bzLG\cegalZ\bzs\text{ or }\bwLG\cegalW\bws \right\} \right.\right)\\
&=& \Prob\left(\bzLG\cegalZ\bzs\text{ or }\bwLG\cegalW\bws \right)  + \Prob\left(\dHinf\left(\permbthetaLG,\bthetas\right)>t, \bzLG\egalZ\bzs, \bwLG\egalW\bws \right)\\
 &\leq& \Prob\left(\bzLG\cegalZ\bzs\right)+\Prob\left(\bwLG\cegalW\bws\right) \\
 &&\quad+ \Rev{\Prob\left(\dHinf\left(\permbthetaLG,\bthetas\right)>t\left|\bzLG\egalZ\bzs,\bwLG\egalW\bws\right.\right)\Prob\left(\bzLG\egalZ\bzs,\bwLG\egalW\bws\right)}\\
& \Rev{\leq}& \Rev{ \Prob\left(\bzLG\cegalZ\bzs\right)+\Prob\left(\bwLG\cegalW\bws\right) + \Prob\left(\dHinf\left(\permbthetaLG,\bthetas\right)>t\left|\bzLG\egalZ\bzs,\bwLG\egalW\bws\right.\right)}\\
\end{eqnarray*}

To complete the proof, we then need to bound from above the terms of this inequality. The two first terms are bounded using Proposition \ref{prop:concentration_inequality:g} in Appendix \ref{Appendix:prop:concentration_inequality:g}, and the last term is bounded with Proposition \ref{prop:btheta} in Appendix \ref{Appendix:prop:btheta}.

\subsection{Concentration inequality on $\bzLG$}\label{Appendix:prop:concentration_inequality:g}

Let us first define the following events.
\begin{itemize}
    \item There is at least one individual in each row class, denoted by
    \[\Ags=\bigcap_{\kk=1}^{\gs}\left\{\Zssumk\neq0\right\}.\]
    \item Denoting $\D$ the random variable equal to the maximal distance between $\Tid$ and the center of the class of row $i$:
    \[\D=\underset{1\leq\kk\leq\gs}{\max}\;\underset{\overset{1\leq\ii\leq\n}{\underset{\text{with }\zsik=1}{}}}{\sup}\;\left|\Tid-\tauk\right| ,\]
    we also define:
    \[\ASg=\{2\D < \Sg < \deltapi-2\D\} \text{ and } \Aid=\Ags \cap \ASg. \] 
\end{itemize}

\begin{lemma}[\Rev{Interesting event}]\label{lemma:event:g}
\[ \Aid \subset \{ \gLG = \gs \} \cap \{ \bzLG \egalZ \bzs \} \]
\end{lemma}

\begin{proof}
On the event $\ASg$, for any two rows $\ii\neq\ip\in\{1,\ldots,\n\}$, we have two possibilities:
\begin{itemize}
    \item Either the rows $\ii$ and $\ip$ are in the same class $\kk$, and then on $\ASg$, we have:
    \[
        \left|\Tid -\Tipd\right| \leq \left|\Tid -\tauk\right|+\left|\Tipd-\tauk\right| \leq 2 \D < \Sg.\\
    \]
    \item Or row $\ii$ is in the class $\kk$ and row $\ip$ in the class $\kp\neq \kk$, and on the event $\ASg$, we have:
    \begin{eqnarray*}
        \left|\Tid -\Tipd\right|&=&\left|\Tid -\taukp-\left(\Tipd-\taukp\right)\right|\\
                                &\geq&\left|\Tid -\taukp\right|-\left|\Tipd-\taukp\right|\\
                                &\geq&\left|\Tid -\taukp\right|-\D\\
                                &\geq&\left|\tauk-\taukp\right|-\left|\Tid -\tauk\right|-\D\\
                                &\geq&\deltapi-2\D\\
                                & > &\Sg.\\
    \end{eqnarray*}
\end{itemize}

Therefore, $G_{\ii}=\Tipard-\Tiparund$ is less than $\Sg$ if and only if both rows $(\ii-1)$ and $(\ii)$ are in the same class. On $\ASg$, the algorithm hence finds the true classification. Moreover, on $\Ags$, there is at least one row in each class, then the algorithm finds the true number of classes. As a conclusion, on $\Aid$, both $\gLG=\gs$ and $\bzLG\egalZ\bzs$ are satisfied.
\end{proof}

\begin{lemma}[\Rev{Concentration inequality for $\Aid$}]\label{lemma:bound:g}
Under Assumption \eqref{hyp:identifiability} and if $\Sg\in ]0,\deltapi[ \text{ and } \Sm \in ]0,\deltarho[$ for $\n,\dd$ large enough:
 \[ \Prob\left(\overline{\Aid}\right) \leq 2 \n \exp\left(- \frac{\dd}{2} \min(\deltapi-\Sg,\Sg)^2 \right) +\gs\left(1-\pimin\right)^{\n} \]
 \Rev{where $\overline{\Aid}$ is the complementary of the event $\Aid$.}
\end{lemma}

\begin{proof}
    Using an union bound, we first obtain:
\[ \Prob\left(\overline{\Aid}\right) \leq \Prob\left(\overline{\Ags}\right) + \Prob\left(\overline{\ASg}\right) \]

Now we bound from above each of these terms. Again with an union bound:

\begin{align*}
\Prob\left(\overline{\Ags}\right)&=\Prob\left(\bigcup_{\kk=1}^{\gs}\overline{\left\{\Zssumk\neq0\right\}}\right) \\ 
				 & \leq \sum_{\kk=1}^{\gs}\Prob\left(\overline{\left\{\Zssumk\neq0\right\}}\right)\\
				 & \leq \sum_{\kk=1}^{\gs}\Prob\left(\Zssumk=0\right)\\
				 &\leq \sum_{\kk=1}^{\gs}\prod_{\ii=1}^{\n} \Prob\left(\Zsik=0\right)\\
				 &\leq \sum_{\kk=1}^{\gs}\prod_{\ii=1}^{\n}\left(1-\piks\right)\\
                                 &\leq\sum_{\kk=1}^{\gs}\prod_{\ii=1}^{\n}\left(1-\pimin\right)\\
                                 &\leq\gs\left(1-\pimin\right)^{\n},\\
\end{align*}

\noindent which gives the upper bound of the first term. Secondly:
\begin{eqnarray*}
\ASg &=& \{ 2\D < \Sg < \deltapi - 2 \D \}\\
&=& \{2\D < \Sg,2\D < \deltapi - \Sg \}\\
& =& \left\{ \D < \frac12 \min(\deltapi-\Sg,\Sg) \right\}.\\
\end{eqnarray*}

Denoting $t= \min(\deltapi -\Sg,\Sg)$,
\begin{align*}
\Prob\left( \overline{\ASg} \right) & = \Prob\left(  \D \geq \frac t2 \right)\\
		       &= \Prob\left(\bigcup_{\kk=1}^{\gs}\bigcup_{\ii|\Rev{\zsik}=1}\left\{\left|\Tid -\tauk\right| \geq \frac t2\right\}\right)\\
                      & \leq\sum_{\kk=1}^{\gs}\sum_{\ii|\Rev{\zsik}=1}\Prob\left(\left|\Tid -\tauk\right| \geq \frac t2 \right).
\end{align*}

Moreover for all $\ii\in\{1,\ldots,\n\}$, given $\zsik=1$, $\Xiplus$ has a binomial distribution $\Bin\left(\dd,\tauk\right)$. The concentration properties of this distribution are then exploited through the Hoeffding inequality:
\begin{eqnarray*}
\Prob\left(\left|\Tid -\tauk\right| \geq \frac t2 \right) = \Prob\left(\left|\Xiplus -\dd\tauk\right| \geq \frac{\dd t}{2} \right) &\leq& 2e^{-\frac12 \dd t^2}.
\end{eqnarray*}

And as a conclusion, the bound of the second term is:
\[ \Prob\left( \overline{\ASg} \right) \leq \sum_{\kk=1}^{\gs}\sum_{\ii|\Rev{\zsik}=1} 2e^{-\frac12 \dd t^2} = 2\n e^{-\frac12 \dd t^2}. \]
\end{proof}

\Rev{With these two lemmas, the following proposition is obtained}
\begin{proposition}\label{prop:concentration_inequality:g}\
Under Assumption \eqref{hyp:identifiability} and if $\Sg\in ]0,\deltapi[ \text{ and } \Sm \in ]0,\deltarho[$ for $\n,\dd$ large enough:
\begin{eqnarray*}
\Prob\left(\bzLG\cegalZ\bzs\right)
&\leq& 2\n \exp\left( -\frac{\dd}{2} \min(\deltapi-\Sg,\Sg)^2 \right) +\gs\left(1-\pimin\right)^\n.\\
\end{eqnarray*}
\begin{eqnarray*}
\Prob\left(\bwLG\cegalW\bws\right)
&\leq& 2\dd \exp\left( -\frac{\n}{2} \min(\deltarho-\Sm,\Sm)^2 \right) +\ms\left(1-\rhomin\right)^\dd.\\
\end{eqnarray*}
\end{proposition}

\begin{proof}

Lemma \ref{lemma:event:g} tells that whenever the event $\Aid$ is satisfied, then both true number of row classes and their true classification are obtained. Lemma \ref{lemma:bound:g} provides an upper bound of $\Prob\left(\overline{\Aid}\right)$. From these lemmas, it is directly deduced that:

\begin{eqnarray*}
\Prob \left(  \bzLG \cegalZ \bzs  \right)
 &\leq& \Prob \left( \overline{\Aid} \right)\\
  &\leq& 2 \n \exp\left(- \frac{\dd}{2} \min(\deltapi-\Sg,\Sg)^2 \right) +\gs\left(1-\pimin\right)^{\n} ,\\
\end{eqnarray*}
which is Proposition \ref{prop:concentration_inequality:g}.
  
\end{proof}

\subsection{\Rev{Concentration inequality on $\dHinf\left(\permbthetaLG,\bthetas\right)>t$}}\label{Appendix:prop:btheta}

\Rev{In this part, the proof of the inequality on $\dHinf\left(\permbthetaLG,\bthetas\right)>t$ is detailled.}

\begin{proposition}\label{prop:btheta}
For all $t>0$, we have:
\begin{eqnarray*}
 &&\!\!\!\!\!\!\!\!\!\!\!\!\Prob\left(\dHinf\left(\permbthetaLG,\bthetas\right)>t\left|\bzLG\egalZ\bzs,\bwLG\egalW\bws\right. \right) \\
 &\leq& \Rev{2\gs\ms\left[1-\pimin\rhomin\left(1-e^{-2t^2}\right)\right]^{nd}}+ 2\gs e^{-2\n t^2}+2\ms e^{-2\dd t^2}\\
\end{eqnarray*}
\end{proposition}

The proof consists in obtaining three bounds: one for each parameter. The inequalities on $\bpi$ and $\brho$ are an application of the Hoeffding inequality and are similar to~\cite{channarond2011Classification} for the row class proportions. To obtain the inequality for $\bal$, it is necessary to study the conditional probability, given the true partition $\left(\bzs,\bws\right)$. Apart from the problem of two asymptotic behaviors, the proof is similar to~\cite{channarond2011Classification}.\\

In the sequel, and for ease of reading, we remove the superscripts $\perms_{\clZ}$ and $\permt_{\clW}$.
Therefore, for all $t>0$:
\begin{eqnarray*}
&&\!\!\!\!\!\!\!\!\!\!\!\!\Prob\left(\dHinf\left(\bthetaLG,\bthetas\right)>t\left|\bzLG\egalZ\bzs,\bwLG\egalW\bws\right.\right)\\
&=&\Prob\left(\max\left(\norminf{\bpiLG-\bpis},\norminf{\brhoLG-\brhos},\norminf{\balLG-\bals}\right)>t \left|\bzLG\egalZ\bzs,\bwLG\egalW\bws\right.\right)\\
&\leq&\Prob\left(\norminf{\bpiLG-\bpis}>t\left|\bzLG\egalZ\bzs,\bwLG\egalW\bws\right.\right)+ \Prob\left(\norminf{\brhoLG-\brhos}>t\left|\bzLG\egalZ\bzs,\bwLG\egalW\bws\right.\right)\\
& &+ \Prob\left(\norminf{\balLG-\bals}>t\left|\bzLG\egalZ\bzs,\bwLG\egalW\bws\right.\right)\\
&\leq&\sum_{\kk=1}^{\gs}\Prob\left(\left|\pikLG-\piks\right|>t\left|\bzLG\egalZ\bzs,\bwLG\egalW\bws\right.\right)+ \sum_{\el=1}^{\ms}\Prob\left(\left|\rholLG-\rhols\right|>t\left|\bzLG\egalZ\bzs,\bwLG\egalW\bws\right.\right)\\
& &+ \sum_{\kk=1}^{\gs}\sum_{\el=1}^{\ms}\Prob\left(\left|\alklLG-\alkls\right|>t\left|\bzLG\egalZ\bzs,\bwLG\egalW\bws\right.\right).
\end{eqnarray*}

The upper bounds of the first and second terms are the same as~\cite{channarond2011Classification}; only the last term is different. For $\alklLG$, first note that when $\bzLG\egalZ\bzs$ and $\bwLG\egalW\bws$
\[\alklLG=\alkltilde=\frac{1}{\zssumk\wssuml}\sum_{\left(\ii,\jj\right)\left|\zsik\wsjl=1\right.}\Xij\]
and given $\left(\bzs,\bws\right)$, the Hoeffding inequality gives for all $t>0$:
\begin{eqnarray*}
\Prob\left(\left|\alklLG-\alkls\right|>t\left|\bzLG\egalZ\bzs,\bwLG\egalW\bws\right.\right)
&=&\Prob\left(\left|\alkltilde-\alskl\right|>t\left|\bzLG\egalZ\bzs,\bwLG\egalW\bws\right.\right)\\
&\leq&\Prob\left(\left|\alkltilde-\alskl\right|>t\right) \\
&\leq&\EspC{\bZs,\bWs}{\left.\Prob\left(\left|\alkltilde-\alskl\right|>t\right|\bZs,\bWs\right)}\\
&\leq&\EspC{\bZs,\bWs}{2e^{-2\Zssumk\Wssuml t^2}}.\\
\end{eqnarray*}
\Rev{But, as $\zssumk=\sum_{i=1}^n\zsik$ and $\wssuml=\sum_{j=1}^d\wsjl$ and the variables are independents, the expectation is:
\begin{eqnarray*}
    \EspC{\bZs,\bWs}{2e^{-2\Zssumk\Wssuml t^2}}&=&2\EspC{\bZs,\bWs}{\exp\left(-2\sum_{i=1}^n\Zsik\sum_{j=1}^d\Wsjl t^2\right)}\\
    &=&\EspC{\bZs,\bWs}{2\prod_{i=1}^n\exp\left(-2\Zsik\sum_{j=1}^d\Wsjl t^2\right)}\\
    &=&2\prod_{i=1}^n\EspC{\bZs,\bWs}{\exp\left(-2\Zsik\sum_{j=1}^d\Wsjl t^2\right)}\\
    &=&2\prod_{i=1}^n\prod_{j=1}^d\EspC{\bZs,\bWs}{\exp\left(-2\Zsik\Wsjl t^2\right)}\\
\end{eqnarray*}
Since, for all $i\in\{1,\ldots,n\}$, $j\in\{1,\ldots,d\}$, $k\in\{1,\ldots,\gs\}$ and $\ell\in\{1,\ldots,\ms\}$, the variable $\Zsik\Wsjl$ is a Bernoulli of parameter $\piks\rhols$, then:
\begin{eqnarray*}
    \EspC{\bZs,\bWs}{\exp\left(-2\zsik\wsjl t^2\right)}&=&\piks\rhols\exp\left(-2\times 1\times t^2\right)+\left(1-\piks\rhols\right)\exp\left(-2\times 0\times t^2\right)\\
    &=&\piks\rhols\exp\left(-2t^2\right)+1-\piks\rhols\\
    &=&1-\piks\rhols\left(1-e^{-2t^2}\right).\\
\end{eqnarray*}
Finally, the inequality is:
\begin{eqnarray*}
    \Prob\left(\left|\alklLG-\alkls\right|>t\left|\bzLG\egalZ\bzs,\bwLG\egalW\bws\right.\right)&\leq&\left[1-\piks\rhols\left(1-e^{-2t^2}\right)\right]^{nd}\\
    &\leq&\left[1-\pimin\rhomin\left(1-e^{-2t^2}\right)\right]^{nd}\\
\end{eqnarray*}
and the result is obtained by the sum on each cluster.
}

\section{Proof of Theorem~\ref{th:consistency}: consistency }\label{cor:Appendix:consistency}

The proof is based on Theorem~\ref{thm:concentration}, as $\n\to+\infty$ and $\dd\to+\infty$ \Rev{and by the assumption~\eqref{hyp:identifiability}, the following limits are obtained for every $t>0$:}
\begin{eqnarray*}
\gs\left(1-\pimin\right)^\n+\ms\left(1-\rhomin\right)^\dd&\underset{\n,\dd\to+\infty}{\longrightarrow}&0\\
\text{and }\quad\gs e^{-2\n t^2}+\ms e^{-2\dd t^2}&\underset{\n,\dd\to+\infty}{\longrightarrow}&0.\\
\end{eqnarray*}
\Rev{For the same reasons, as soon as $t$ is positive, $0<\pimin\rhomin\left(1-e^{-2t^2}\right)<1$ and
\begin{eqnarray*}
\gs\ms\left[1-\pimin\rhomin\left(1-e^{-2t^2}\right)\right]^{nd}
&\underset{\n,\dd\to+\infty}{\longrightarrow}&0.\\
\end{eqnarray*}}
For the last terms, the assumption~\eqref{hyp:vanishing} gives that
\Rev{
\begin{eqnarray*}
    \frac{\n}{\log(\dd)}\underset{\n,\dd\to+\infty}{\longrightarrow}+\infty&\text{ and }&\frac{\dd}{\log(\n)}\underset{\n,\dd\to+\infty}{\longrightarrow}+\infty\\
\end{eqnarray*}
which results in}
\begin{eqnarray*}
\n e^{-\frac{1}{2}\dd\deltapi^2}+\dd e^{-\frac{1}{2}\n\deltarho^2}&\underset{\n,\dd\to+\infty}{\longrightarrow}&0\\
\end{eqnarray*}
\Rev{and, as $\Sgnd$ tends to zero,}
\[\n e^{-\dd \frac{\left(\deltapi -\Sgnd\right)^2}{2}}\underset{\n,\dd\to+\infty}{\longrightarrow}0.\]
\Rev{Thanks the same Assumption \eqref{hyp:vanishing},} there exists a positive constant $C>\sqrt{2}$ such that for $\n$ and $\dd$ large enough
\[\Sgnd\sqrt{\frac{\dd}{\log\n}}>C\Longrightarrow\frac{\Sgnd}{\sqrt{2}}\sqrt{\frac{\dd}{\log\n}}>\frac{C}{\sqrt{2}}>1\]

\begin{eqnarray*}
\n e^{-\dd \frac{{\Sgnd}^2}{2}}&=&\exp\left[\log\n-\dd \frac{{\Sgnd}^2}{2}\right]\\
                            &=&\exp\left[\log\n\left(1-\left(\sqrt{\frac{\dd}{\log\n}}\frac{\Sgnd}{\sqrt{2}}\right)^2\right)\right]\\
                            &\leq&\exp\left[\log\n\underbrace{\left(1-\frac{C}{\sqrt{2}}\right)}_{<0}\right]\\
                            &\underset{\n,\dd\to+\infty}{\longrightarrow}&0.\\
\end{eqnarray*}

\section{\Rev{Proof of Theorem~\ref{th:consistency:spare}: consistency in the sparse case}}\label{cor:Appendix:consistency:spare}

\Rev{First, if $\piminnd$ tends to zero (and as $\gsnd\leq 1/\piminnd$), the series expansion of $t\mapsto \log(1-t)$ is used
\begin{eqnarray*}
    \gsnd\left(1-\piminnd\right)^\n&=&\exp\left[\log\gsnd+\n\log\left(1-\piminnd\right)\right]\\
    &\leq&\exp\left\{\n\left[\frac{1}{\n}\log\left(\frac{1}{\piminnd}\right)+\log\left(1-\piminnd\right)\right]\right\}\\
    &\leq&\exp\left\{\n\left[-\frac{1}{\n}\log\left(\piminnd\right)-\piminnd+o\left(\piminnd\right)\right]\right\}\\
    &\leq&\exp\left\{\n\left[-\piminnd+o\left(\piminnd\right)\right]\right\}\\
    &\leq&\exp\left[-\n\piminnd+o\left(\n\piminnd\right)\right]\\
\end{eqnarray*}
and, by the assumption~\eqref{hyp:varying2}, $\n\piminnd$ tends to infinity and 
\[\gsnd\left(1-\piminnd\right)^\n\underset{\n,\dd\to+\infty}{\longrightarrow}0.\]
For the same reasons, for all $t>0$
\begin{eqnarray*}
    \msnd\left(1-\rhominnd\right)^\n&\underset{\n,\dd\to+\infty}{\longrightarrow}&0\\
    \text{and }\gs\ms\left[1-\pimin\rhomin\left(1-e^{-2t^2}\right)\right]^{nd}&\underset{\n,\dd\to+\infty}{\longrightarrow}&0.\\
\end{eqnarray*}
Moreover, as
\begin{eqnarray*}
    \gsnd e^{-2\n t^2}&=&\exp\left(\log \gsnd-2\n t^2\right)\\
    &\leq&\exp\left\{-\n\left[\frac{1}{\n}\log\left(\piminnd\right)+2 t^2\right]\right\}\\
    &\leq&\exp\left[-\n2 t^2+o\left(\n\right)\right].\\
\end{eqnarray*}
then
\[\gsnd e^{-2\n t^2}\underset{\n,\dd\to+\infty}{\longrightarrow}0\quad\text{ and }\quad\msnd e^{-2\dd t^2}\underset{\n,\dd\to+\infty}{\longrightarrow}0.\]
Finally, the assumption~\eqref{hyp:varying1} implies that there exists a positive constant $C>2$ such that for $\n$ and $\dd$ large enough
\begin{eqnarray*}
    \frac{\deltapind}{\Sgnd}>C&\Leftrightarrow&\deltapind>\Sgnd C\\
    &\Leftrightarrow&\deltapind>\Sgnd (C-1)+\Sgnd\\
    &\Leftrightarrow&\deltapind-\Sgnd>\Sgnd \underbrace{(C-1)}_{>1}>\Sgnd\\
\end{eqnarray*}
and, for $\n$ and $ \dd$ large enough, $\min\left(\deltapind-\Sgnd,\Sgnd\right)$ is $\Sgnd$ and the assumptions~\eqref{hyp:identifiability} and~\eqref{hyp:vanishing} allow to conclude.
}

\section{Supplementary Material}

In this supplementary material, additional figures from the experiments in the section~\ref{sec:simulation} are presented.

\begin{figure}[!h]
\begin{center}
    \includegraphics[width=\linewidth]{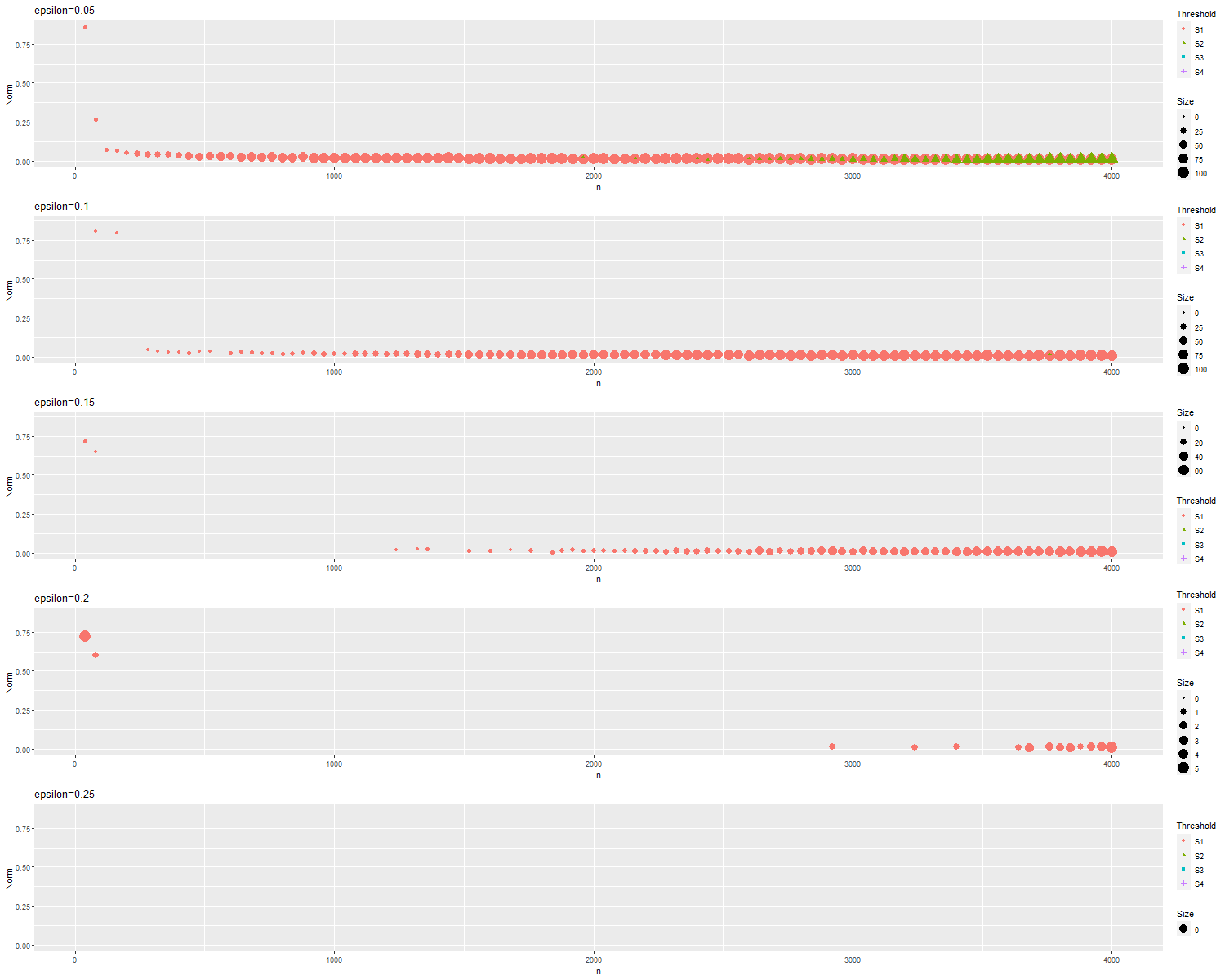}
\caption{\label{Fig:Res:Theta:Neq}Average estimate of the distance $\dHinf\left(\bthetas,\bthetaLG\right)$ between the true parameters and the estimated parameters following $\varepsilon$ (rows) in the arithmetic case: for each graphic, the threshold is represented by different symbols (\textcolor{RubineRed}{$\bullet$} for $\Sun$, \textcolor{LimeGreen}{$\blacktriangle$} for $\Sdeux$, \textcolor{Cerulean}{$\blacksquare$} for $\Strois$ and  \textcolor{purple}{$\boldsymbol{+}$} for $\Strois$) and the size for the number of finite values used; the number of rows $\n$ varies between 40 and 4000 and $\dd$ is supposed equal $\n$.}
\end{center}
\end{figure}

\begin{figure}[!h]
\begin{center}
    \includegraphics[width=\linewidth]{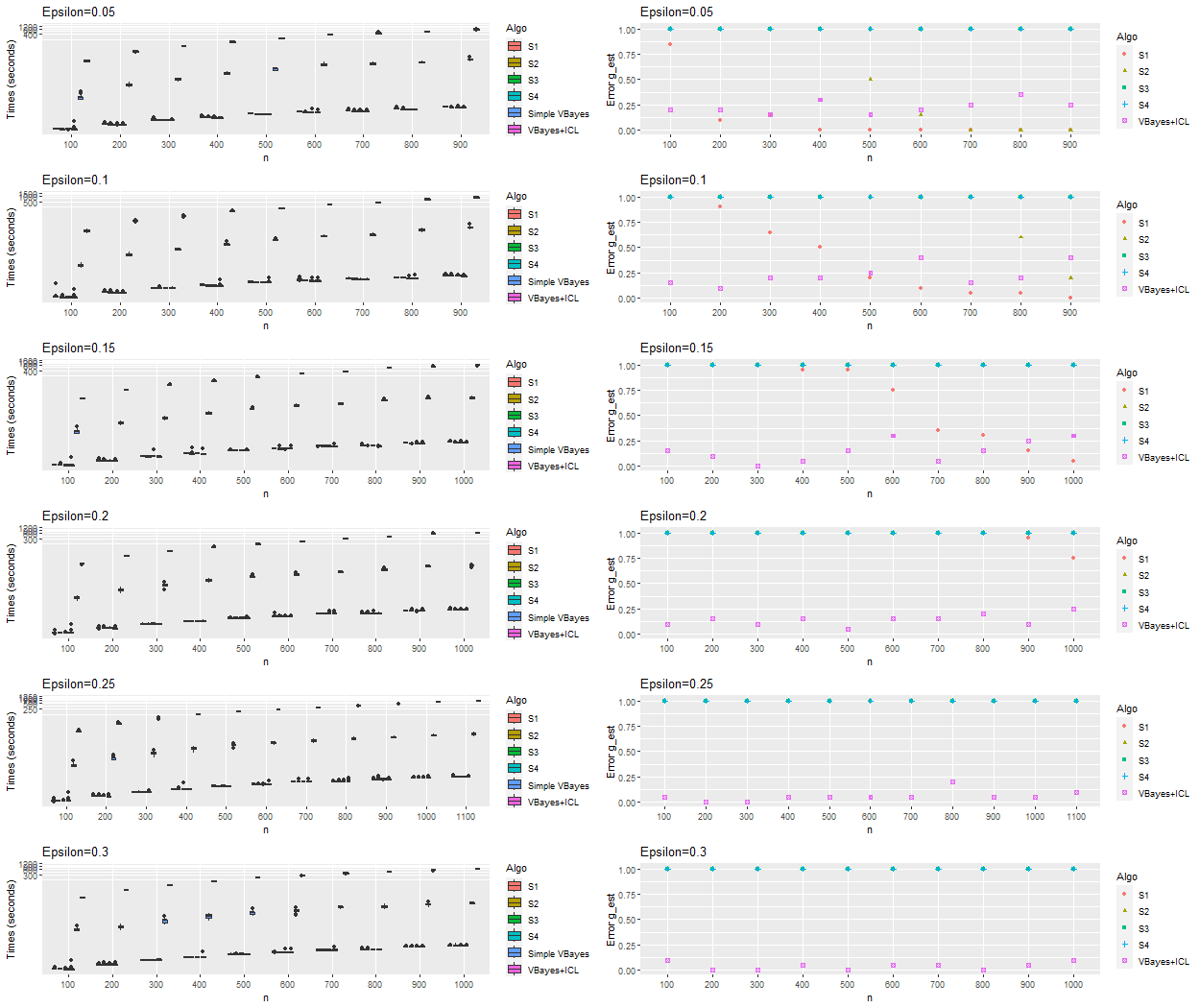}
\caption{\label{Fig:times:Eq}\Rev{On the left, boxplots of the computing times (in seconds; logarithmic scale) for each procedure (colour) in function of the numbers of rows and columns ($\n=\dd$) for each $\varepsilon$ (in rows) over the 200 simulations for the balanced case. On the bottom, boxplots of quality of estimations averaged over the 20 matrices for each $\varepsilon$ (in rows) for each procedure (colour) in function of the numbers of rows and columns ($\n=\dd$) for the balanced case.}}
\end{center}
\end{figure}

\begin{figure}[!h]
\begin{center}
    \includegraphics[width=\linewidth]{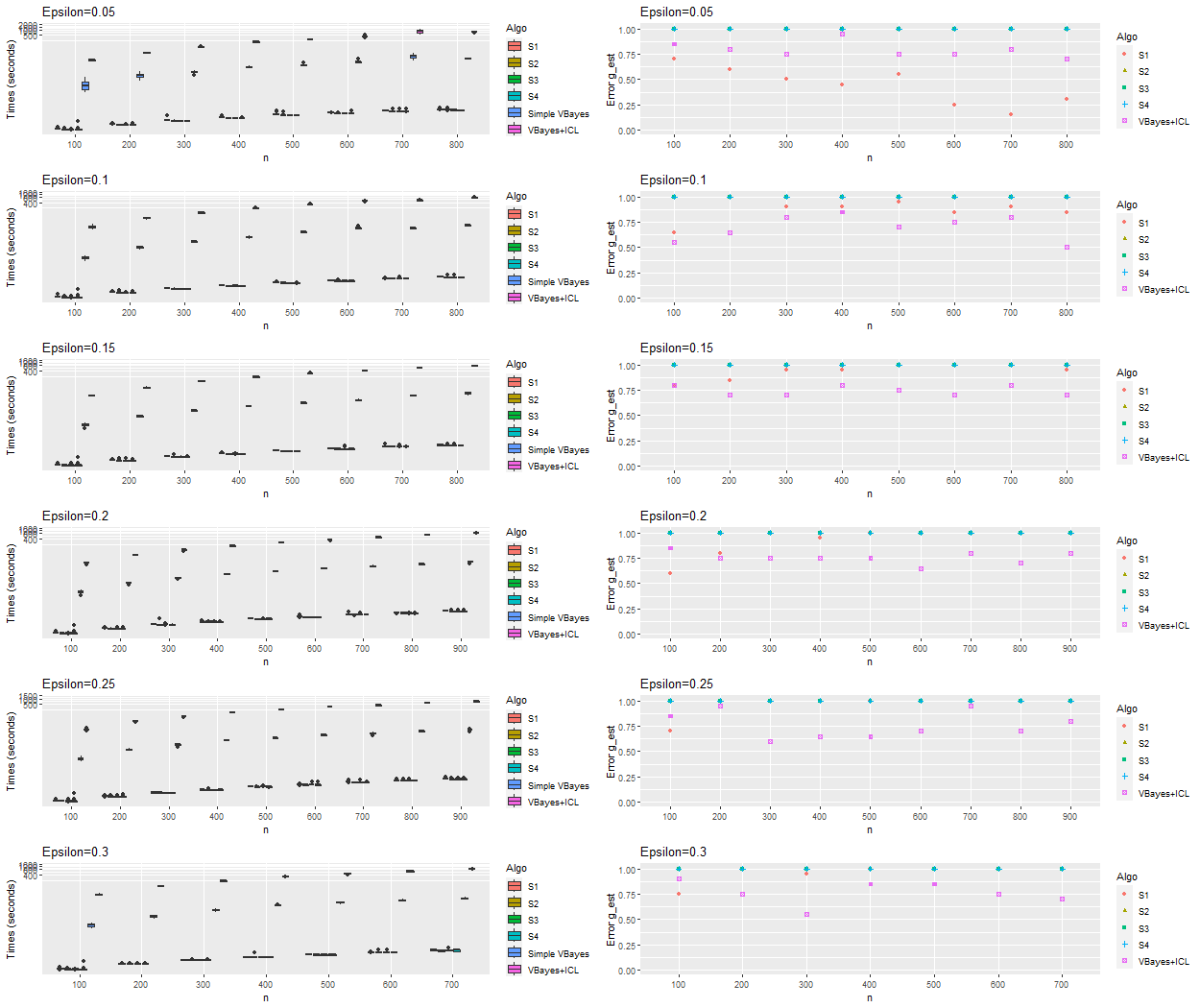}
\caption{\label{Fig:times:Neq}\Rev{On the left, boxplots of the computing times (in seconds; logarithmic scale) for each procedure (colour) in function of the numbers of rows and columns ($\n=\dd$) for each $\varepsilon$ (in rows) over the 200 simulations for the arithmetic case. On the bottom, boxplots of quality of estimations averaged over the 20 matrices for each $\varepsilon$ (in rows) for each procedure (colour) in function of the numbers of rows and columns ($\n=\dd$) for the arithmetic case.}}
\end{center}
\end{figure}


%
%

\bibliographystyle{spmpsci}      
\bibliography{biblio}   

\end{document}